\documentclass[11pt]{amsart}
\usepackage{amsmath,amssymb,delarray,amsthm,bbm}

\usepackage[hmargin=0.95in,height=8.6in]{geometry}

\usepackage{ifpdf}
\ifpdf
\usepackage[pdftex]{graphicx}
\else
\usepackage[dvips]{graphicx}
\fi

\usepackage{ifpdf}
\usepackage{color}
\definecolor{webgreen}{rgb}{0,.5,0}
\definecolor{webbrown}{rgb}{.8,0,0}
\definecolor{emphcolor}{rgb}{0.95,0.95,0.95}

\usepackage{hyperref}
\hypersetup{%
          colorlinks=true,
          linkcolor=webbrown,
          filecolor=webbrown,
          citecolor=webgreen,
          breaklinks=true}
\ifpdf
\hypersetup{pdftex,
            pdfstartview=FitH, 
            bookmarksopen=true,
            bookmarksnumbered=true
}
\else
\hypersetup{dvips}
\fi

\newcommand\E{\ensuremath{\mathbb{E}}}
\newcommand\R{\ensuremath{\mathbb{R}}}
\newcommand\PP{\ensuremath{\mathbb{P}}}
\newcommand\tP{\ensuremath{\tilde{\mathbb{P}}}}
\newcommand\tE{\ensuremath{\tilde{\mathbb{E}}}}

\newcommand\e{\ensuremath{\mathrm{e}}}

\newcommand\defn{\triangleq}   

\newcommand\half{\frac{1}{2}}

\newcommand\Fy{\mathcal{F}^Y}
\newcommand\F{\mathcal{F}}
\newcommand\bW{\overline{W}}

\newcommand\Dt{\Delta t}
\newcommand\tDt{{{t+\Delta t}}}

\newcommand\tT{\tilde{T}}

\newcommand\tFy{\tilde{\mathcal{F}}^Y}
\newcommand\cC{\mathcal{C}}
\newcommand\scE{\mathcal{E}}
\newcommand\cP{\mathcal{P}}
\newcommand\N{\mathcal{N}}
\newcommand\cB{\mathcal{B}}
\renewcommand\H{\mathcal{H}}
\newcommand\cO{\mathcal{O}}
\newcommand\Z{\hat{\mathcal{Z}}}
\newcommand\bZ{\mathbf{Z}}
\newcommand\cZ{\mathcal{Z}}
\newcommand\hq{\hat{q}}
\newcommand\vth{\vartheta}
\DeclareMathOperator*{\pr}{pr}
\DeclareMathOperator*{\tr}{tr}
\DeclareMathOperator*{\argmin}{arg\, min}
\newcommand\prH{\textstyle\pr_{\H}}

\newtheorem{thm}{Theorem}
\newtheorem{lemma}{Lemma}
\newtheorem{prop}{Proposition}
\newtheorem{assume}{Assumption}

\theoremstyle{remark}
\newtheorem{rem}{Remark}

\begin{document}
\title{A Simulation Approach to Optimal Stopping Under Partial Information}
\author{Mike Ludkovski}
\address{Department of Statistics and Applied Probability\\
University of California\\
Santa Barbara, CA 93106 USA \\ Phone: 1(805)893-5634. Email: \url{ludkovski@pstat.ucsb.edu}}

\date{\today}
\begin{abstract}
We study the numerical solution of nonlinear partially observed optimal stopping problems. The
system state is taken to be a multi-dimensional diffusion and drives the drift of the
observation process, which is another multi-dimensional diffusion with correlated noise.
Such models where the controller is not fully aware of her environment are of interest in applied probability and financial mathematics. We propose a new approximate numerical algorithm based on the particle filtering and regression Monte Carlo methods. The algorithm maintains a continuous state-space and yields an integrated approach to the filtering and control sub-problems. Our approach is entirely simulation-based and therefore allows for a robust implementation with respect to
model specification. We carry out the error analysis of our scheme and illustrate with several
computational examples. An extension to discretely observed stochastic volatility models is also considered.
\end{abstract}

\keywords{optimal stopping, nonlinear filtering, particle filter, Snell envelope, regression Monte Carlo}

\maketitle

\section{Introduction}

Let $(\Omega,\mathcal{F}, (\F_t), \PP)$ be a filtered probability space and consider a
$d$-dimensional process $X = (X_t)$ satisfying an \^{I}to stochastic differential equation (SDE) of the form
\begin{align}\label{eq:X-sde}
dX_t = b(X_t) \, dt + \alpha(X_t) \, dU_t + \sigma(X_t) \, dW_t,
 \end{align}
where $U$ and $W$ are two independent $(\F_t)$-adapted Wiener processes of dimension $d_U$ and $d_W$ respectively.
Let $Y$ be a $d_Y \equiv d_U$-dimensional diffusion given by
\begin{align}\label{eq:Y-sde}
dY_t = h(X_t) \, dt + \, dU_t.
\end{align}
Assumptions about the coefficients of \eqref{eq:X-sde}-\eqref{eq:Y-sde} will be given later.
Denote by $\Fy_t = \sigma( Y_s \colon 0 \le s \le t)$ the filtration generated by $Y$. We study
the partially observed finite horizon optimal stopping problem
\begin{align}\label{eq:basic-problem}
\sup_{\tau \le T, \; \Fy-\text{adapted}} \E \left[ g(\tau, X_\tau, Y_\tau) \right],
\end{align}
where $g: [0,T] \times \R^{d} \times \R^{d_Y} \to \R$ is the reward functional.

The probabilistic interpretation of \eqref{eq:basic-problem} is as follows. A controller wishes to maximize expected reward $g(t,x,y)$ by selecting an optimal stopping time $\tau$. Unfortunately, she only has access to the observation process $Y$; the state $X$ is not revealed and can be only partially inferred through its impact on the drift of $Y$. Thus, $\tau$ must be based on the information contained solely in $Y$. Recall that even when $Y$ is observed continuously, its drift is never known with certainty; in contrast the instantaneous volatility of $Y$ can be obtained from the corresponding quadratic variation.

Such partially observed problems arise frequently in financial
mathematics and applied probability where the agent is not fully aware of her
environment, see Section \ref{sec:apps} below. One of their interesting features is the interaction between learning and
optimization. Namely, the observation process $Y$ plays a \emph{dual}
role as a source of information about the system state $X$, and as a
reward ingredient.  Consequently, the agent has to consider the trade-off between further monitoring of $Y$ in order to obtain a more accurate inference of $X$, vis-a-vis stopping early in case the state of the world is unfavorable. This tension between exploration and maximization is even more accentuated when time-discounting is present. Compared to the fully observed setting, we therefore expect that partial information would postpone decisions due to the demand for learning.

In the given form the problem \eqref{eq:basic-problem} is non-standard because the payoff $g(t,X_t,Y_t)$ is not adapted to the observed filtration $(\Fy_t)$ and, moreover, $Y$ is not Markovian with respect to
$(\Fy_t)$. This difficulty is resolved by a two-step inference/optimization approach. Namely, the first filtering step transforms \eqref{eq:basic-problem} into an equivalent fully-observed formulation
using the Markov conditional distribution $\pi_t$ of $X_t$ given $\Fy_t$. In the second step, the resulting standard optimal stopping problem with the Markovian state $(\pi_t, Y_t)$ is solved.

Each of the two sub-problems above are covered by an extensive literature. The filtering problem with diffusion observations was first studied by Kalman and Bucy \cite{KalmanBucy} and we refer to the excellent texts \cite{BensoussanBook,KallianpurBook} for the general theory of nonlinear stochastic filtering. The original linear model of \cite{KalmanBucy} had a key advantage in the availability of sufficient statistics and subsequent closed-form filtering formulas for $\pi_t$. Other special cases where the filter was explicitly computable were obtained by \cite{Benes} and \cite{BrigoHanzonLegland99}. However, in the general setup of \eqref{eq:X-sde}-\eqref{eq:Y-sde}, the conditional distribution $\pi_t$ of $X_t$ is measure-valued, i.e.\ an infinite-dimensional object. This precludes consideration of explicit solutions and poses severe computational challenges.

To address such nonlinear models, a variety of approximation tools have been proposed. First, one may linearize the system \eqref{eq:X-sde}-\eqref{eq:Y-sde} by applying (A) the \emph{extended Kalman filter} \cite{ItoXiong00,BudhirajaKushnerIEEE}. Thus, the conditional distribution of $X$ is summarized by its conditional mean $m_t = \E[X_t | \Fy_t]$ and conditional variance $P_t = \E[ (X_t - m_t)^2 | \Fy_t]$. One then derives (approximate) evolution equations for $(m_t, P_t)$ given observations $Y$. More generally, $\pi_t$ can be parameterized by a given family of probability densities, yielding the (B) \emph{projection filter}. Let us especially single out the exponential projection methods studied by Brigo et al.\ \cite{BrigoHanzonLegland99,BrigoHanzonIEEE98}. Third, the state space of $\pi_t$ can be discretized through (C) \emph{optimal quantization} methods \cite{PagesPham,SellamiPhamRunggaldier}. This replaces $(\pi_t)$ by a non-Markovian approximation $(\tilde{\pi}_t)$ whose transition probabilities are pre-processed via Monte Carlo simulation.  Fourth, one may apply (D) \emph{Wiener chaos expansion} methods \cite{Lototsky06,LototskyMikulevicius97,RozovskiiMikulevicius99} that reduce computation of $\pi_t$ to a solution of SDE's plus ordinary differential equation systems. Finally, (E) \emph{interacting particle systems} have been considered to approximate $\pi_t$ non-parametrically via simulation tools \cite{Crisan06,Gaines98,CrisanLyons99,DelMoralProtter01}.

The optimal stopping sub-problem of the second step can again be tackled within several  frameworks. When the transition density of the state variables is known, classical (a) \emph{dynamic programming} computations are possible, see e.g.\ \cite{Sh78}. If the problem state is low-dimensional and Markov, one may alternatively use the quasi-variational formulation to obtain a free-boundary partial differential equation (pde) and then implement a (b) \emph{numerical pde solver} for an efficient solution. Thirdly, (c) \emph{simulation-based methods} \cite{Egloff05,Longstaff,TsitsiklisVanRoy00} that rely on probabilistic Snell envelope techniques can be applied.

The joint problem of optimal stopping with partial observations was treated in \cite{GatarekSwiech}, \cite{GozziSwiechRouy}, \cite{Mazziotto86}, \cite{SellamiPhamRunggaldier} and \cite{SmithMoscarini}. All these models can be viewed as a combination of the listed approaches to the two filtering/optimization sub-problems. For example, \cite{Mazziotto86} proposes to use the assumed density filter for the filtering step, followed by a pde solver for the optimization. This can be summarized as algorithm (B)/(b) in our notation. Meanwhile, \cite{SellamiPhamRunggaldier} use (C)/(a), i.e.\ optimal quantization for the filter and then dynamic programming to find optimal stopping times. Methodologically, two ideas have been studied. First, using
filtering techniques (A) or (B), one may replace $\pi_t$ by a low-dimensional Markovian
approximation $\tilde{\pi}_t$. Depending on the complexity of the model, algorithms (a) or (b)
can then be applied in the second step. Unfortunately, the resulting filtering equations are
inconsistent with the true dynamics of $\pi_t$, and require a lot of computations to
derive them for each considered model. The other alternative is to use the quantization technique (C) which is robust and produces a consistent (but non-Markovian) approximation $\tilde{\pi_t}$. Since the state space of $\tilde{\pi_t}$  is fully discretized, the resulting optimal stopping problem can be solved exactly using dynamic programming algorithm (a). Moreover, tight error bounds are available. The shortcomings of this approach are the need to discretize the state space of $X$ and the requirement of offline pre-processing to compute the transition
density of $\tilde{\pi}_t$.

In this paper we propose a new approach of type (E)/(c) that uses a \emph{particle filter} for the inference step and a \emph{simulation tool} in the optimization step. Our method is attractive based on three accounts. Firstly, being entirely simulation-based it can be generically applied to a wide variety of models, with only minor modifications. In particular, the implementation is robust and requires only the ability to simulate $(X_t, Y_t)$. For comparison, free boundary pde solvers of type (b) often use advanced numerical techniques for stability and accuracy purposes and must be re-programmed for each class of models. Also, in contrast to optimal quantization, no pre-processing is needed. Moreover, the interacting particle system approach to filtering is also robust with respect to different observation schemes. In the original system \eqref{eq:X-sde}-\eqref{eq:Y-sde} it is assumed that $Y$ is observed continuously. It is straightforward to switch our algorithm to discrete regularly-spaced observations of $Y$ that may be more natural in some contexts.

Secondly, our approach maintains a continuous state space throughout all computations. In
particular, the computed optimal stopping rule $\tau^*$ is continuous, eliminating that source of error and leading to a more natural decision criteria for the controller. Thus, compared to optimal quantization, our approach is expected to produce more ``smooth'' optimal stopping boundaries.
Third, our method allows the user to utilize her \emph{domain knowledge} during the optimization step. In most practical applications, the user already has
a guess regarding an optimal stopping rule and the numerical computations are used as a
refinement and precision tool. However, most optimal stopping algorithms rely on a ``brute force'' scheme to obtain an optimal stopping rule. By permitting custom input for the optimization step, our scheme should heuristically lead to reduced computational efforts and increased accuracy.

Finally, maintaining the simulation
paradigm throughout the solution allows us to integrate the filtering and Snell envelope
computations. In particular, by carrying along a high-dimensional approximation of $\pi_t$, the initial filtering errors can be minimized in a flexible and anticipative way with respect to the subsequent optimization step. Thus, the introduction of filtering errors is delayed for as long as possible.  This is important for optimal stopping where the forward-propagated errors (such as the filtering error) strongly affect the subsequent backward recursion solution for $\tau^*$.
To summarize, our scheme should be viewed as an even more flexible alternative for the optimal quantization method of \cite{SellamiPhamRunggaldier}. 

\begin{rem}
To our knowledge the idea of integrated stochastic filtering and optimization was
conceived in \cite{MullerEtal04}, in the context of utility maximization with partially
observed state variables. Muller et al.~\cite{MullerEtal04} proposed to use the Markov Chain Monte Carlo (MCMC) methods and an auxiliary randomized pseudo-control variable to do both steps at once.
These ideas were then further analyzed in \cite{BrandtStroud,Viens03} for a portfolio
optimization problem with unobserved drift parameter and unobserved stochastic volatility, respectively. In fact, Viens et al.\ \cite{Viens03}
utilized a particle filter but then relied on discretizing the control and observation processes  to obtain a finite-dimensional problem with discrete scenarios. While of the same flavor, this approach must be modified for optimal stopping problems like \eqref{eq:basic-problem}, as the control variable $\tau$ is infinite-dimensional. Indeed, stopping rules $\tau$ are in one-to-one correspondence with stopping regions, i.e.\ subsets of the space-time state space. Such objects do not admit easy discretization. Moreover, the explicit presence of time-dimension as part of our control makes MCMC simulation difficult. Thus, we maintain the probabilistic backward recursion solution method instead.
\end{rem}

The rest of the paper is organized as follows. In Section \ref{sec:filter} we recall the
general filtering paradigm for our model and the Snell envelope formulation of the optimal stopping problem \eqref{eq:basic-problem}. Section \ref{sec:algorithm} describes in detail the new algorithm, including the variance-minimizing
branching particle filter in Section \ref{sec:particle}, and the regression Monte Carlo approach to compute the Snell envelope in Section \ref{sec:rms}. We devote Section \ref{sec:errors} to the error analysis of our scheme and to the proof of the overall convergence of the algorithm.
Section \ref{sec:examples}  then illustrates our scheme on a numerical example; a further computational example is provided in Section \ref{sec:stoch-vol} which discusses the extension of our method to discretely observed stochastic volatility models. Finally, Section \ref{sec:conclusion} concludes.

Before proceeding, we now give a small list of applications of the model \eqref{eq:X-sde}-\eqref{eq:basic-problem}.


\subsection{Applications}\label{sec:apps}
%
%

\subsubsection*{Optimal Investment under Partial Information} The following investment timing problem arises in the theory of \emph{real options}. A manager is planning to launch a new project, whose value $(Y_t)$ evolves according to
$$ dY_t = X_t \, dt + \sigma_Y \, dU_t,
$$ where the drift parameter $(X_t)$ is unobserved and $(U_t)$ is an $\R$-valued Wiener process. The environment variable $X_t$ represents the current economic conditions; thus when economy is booming, potential project value grows quickly, whereas it may be declining during a recession. At launch time $\tau$ the received profit $g(\cdot)$ is a function of current project value $Y_\tau$, as well as extra uncertainty that depends on the environment state. For instance, consider $g(\cdot) = Y_\tau \cdot (a_0 + a_1 X_\tau + b_0 \epsilon)$,  $\epsilon \sim \N(0,1)$ independent, where the second term models the profit multiplier based on economy state. Conditioning on the realization of $\epsilon$, expected profit is $g(\tau, X_\tau, Y_\tau) = Y_\tau (a_0 + a_1 X_\tau)$. Such a model with continuous-time observations was considered by \cite{DecampsEtal} in the \emph{static} case where $X_0 \in \{ 0, 1\}$ and $dX_t = 0$. A similar problem was studied in \cite{MiaoWang} with an additional consumption control. 

Using the methods below, we can treat this problem for general $X$-dynamics of the type \eqref{eq:X-sde}, under both continuous and discrete observations.

\subsubsection*{Stochastic Convenience Yield Models}
Compared to holding of financial futures, physical ownership of commodities entails additional benefits and costs. Accordingly,  the rate of return  on the commodity spot contract will be different from the risk-free rate. The stochastic convenience yield models \cite{CL-utah,Schwartz97} postulate
that the drift of the asset price $(Y_t)$ under the pricing measure $\PP$ is itself a stochastic process,
$$ \left\{\begin{aligned} dY_t & = Y_t( X_t \, dt + \sigma_Y \, dU_t), \\
dX_t &= b(X_t) \, dt + \alpha(X_t) \, dU_t + \sigma_X(X_t) \, dW_t. \end{aligned}\right.$$
One may now consider the pricing of American Put options on asset $Y$ with maturity $T$ and strike $K$,
$$ \sup_{\tau \le T} \E[ \e^{-r \tau} (K- Y_\tau)_+ ],$$
where the convenience yield $X$ is unobserved and must be dynamically inferred. Beyond using $Y$ to learn about $X_t$, it is also possible to filter other observables, e.g.\ futures contracts, see \cite{CL-utah}.

\subsubsection*{Reliability Models with Continuous Review}
Quality control models in industrial engineering \cite{JensenHsu} can also be viewed as examples of \eqref{eq:basic-problem}. Let $X_t$ represent the current quality of the manufacturing process. This quality fluctuates due to machinery state and also external disturbances,  such as current workforce effort, random shocks, etc. When quality is high, the revenue stream $Y$ is increasing; conversely poor quality may decrease revenues. Because revenues are also subject to random disturbances, current quality is never observed directly. In this context, it is asked to find an optimal time $\tau$ to replace the machinery (at cost $g(X_t)$) and reset the quality process $(X_t)$. Assuming ``white noise'' shocks to the system and continuous monitoring of revenue stream this leads again to \eqref{eq:X-sde}-\eqref{eq:Y-sde}-\eqref{eq:basic-problem}. The case where $Y$ is discretely observed and $X$ is a finite-state Markov chain was treated by Jensen and Hsu \cite{JensenHsu}. 







\section{Optimization Problem}\label{sec:filter}
\subsection{Notation}
We will use the following notation throughout the paper:
\begin{itemize}
\item For $x \in \R$, we write $x=\lfloor x \rfloor + \{ x \}$ to denote the smallest integer small than $x$ and the fractional part of $x$, respectively.


\item $\delta_{x}$ denotes the Dirac measure at point $x$.

\item $\cC^k_b$ denotes the space of all real-valued, bounded,  continuous functions with bounded continuous  derivatives up to order $k$ on $\R^{d}$. We endow $\cC^k_b(\R^{d})$ with the following norm
$$
\| f \|_{m,\infty} = \sum_{|\alpha| \le m} \sup_{x \in \R^{d}} | D_\alpha f(x)|, \qquad f \in  \cC^k_b(\R^{d}),\qquad m \le k,$$
where $\alpha = (\alpha_1, \ldots, \alpha_{d})$ is a multi-index and  derivatives are written as $D_\alpha f \triangleq \partial_1^{\alpha_1} \cdots \partial_{d}^{\alpha_{d}} f$.

\item $W^k_p = \{ f \colon D_\alpha f \in L^p(\R^d), |\alpha| \le k \}$ denotes the Sobolev space of functions with $p$-integrable derivatives up to order $k$.

\item $\cP(\R^{d})$ is the space of all probability measures over the Borel $\sigma$-algebra $\cB(\R^{d})$. For $\mu \in \cP$, $\mu(f) \triangleq \int_{\R^{d}} f(x) \mu(dx)$. We endow $\cP$ with the weak topology; $\mu_n \to \mu$ weakly if $\forall f \in \cC^0_b$, $\mu_n(f) \to \mu(f)$.

\end{itemize}

\subsection{Filtering Model}
In this section we briefly review the theory of nonlinear filtering as applied to problem \eqref{eq:basic-problem}. We follow \cite{Crisan06} in our presentation.

Before we begin, we make the following technical assumption regarding the coefficients in \eqref{eq:X-sde} and \eqref{eq:Y-sde}.
\begin{assume}
The coefficients of \eqref{eq:X-sde} satisfy: $b(x)\in \cC^3_b(\R^d), \alpha(x) \in \cC^3_b(\R^{d \times d_Y}),\sigma(x) \in \cC^3_b(\R^{d \times d_W}) $ and moreover, $\alpha$ and $\sigma$ are strictly positive-definite matrices of size $d \times d_Y$ and $d \times d_W$ respectively. Similarly, in \eqref{eq:Y-sde}, $h(x) \in \cC^4_b(\R^d)$.
\end{assume}
This assumption in particular guarantees the existence of a unique strong solution to \eqref{eq:X-sde}, \eqref{eq:Y-sde}. We also assume that
\begin{assume}
The payoff function $g$ is bounded and twice jointly continuously differentiable $g \in \cC^2_b([0,T] \times \R^d \times \R^{d_Y})$.
\end{assume}
The latter condition is often violated in practice where payoffs can be unbounded. However, one may always truncate $g$ at some high level $\bar{N}$ without violating the applicability of the model.




We begin by considering the conditional distribution of $X$ given $\Fy_t$. Namely, for $f\in \cC^2_b(\R^{d})$ define
\begin{align}\label{def:pi}
\pi_t f \defn \E[ f(X_t) | \Fy_t].
\end{align}
It is well-known \cite{BensoussanBook} that $\pi_t f$ is a Markov, $\Fy$-adapted process that solves the Kushner-Stratonovich equation
\begin{align}\label{eq:kushner-stratonovich}
d(\pi_t f) = \pi_t( A f) \, dt + \sum_{k=1}^{d_Y} \left[ \pi_t (h_k \cdot f) - \pi_t(h_k) \cdot \pi_t (f) + \pi_t(B^k f) \right] [ dY^k_t - \pi_t(h_k)\, dt],
\end{align}
where the action of the differential operators $A$ and $B^k$ on a test function $f \in \cC^2_b(\R^{d})$ is defined by 
\begin{align}\label{eq:A-B-operators}\left\{ \begin{aligned}
A f(x) & \triangleq \half \sum_{i,j}^{d} \left( \sum_{k=1}^{d_U} \alpha_{ik}(x)\alpha_{jk}(x) + \sum_{k=1}^{d_W}\sigma_{ik}(x)\sigma_{jk}(x) \right)\partial_i \partial_j f(x) + \sum_{i=1}^{d} b_i(x) \partial_i f(x), \\
B^k f(x) & \triangleq \sum_{i=1}^{d} \alpha_{ik}(x)\partial_i f(x).\end{aligned}\right.
\end{align}
Thus, $\pi_t$ is a probability measure-valued process solving the stochastic partial differential equation (spde) corresponding to the adjoint of \eqref{eq:kushner-stratonovich}. To avoid the nonlinearities in \eqref{eq:kushner-stratonovich}, a simpler linear version is obtained by utilizing the reference probability measure device. Define a $\PP$-equivalent probability measure $\tP$ by 
\begin{align}\label{eq:tilde-P}
\frac{ d\tP}{d \PP} \Big|_{\F_T} = \zeta_t \triangleq \exp \left( - \sum_{k=1}^{d_U} \int_0^T h_k(X_s) \, dU^k_s -
\half \sum_{k=1}^{d_U}  \int_0^T h_k^2(X_s) \, ds \right).
\end{align}
From the Girsanov change of measure theorem (recall that $h$ is bounded so that $\E[ \zeta_t] = 1$), it follows that under $\tP$ the observation $Y$ is
a Brownian motion and the signal $X$ satisfies
\begin{align}\label{eq:X-tP}
dX_t = (b(X_t) - \alpha h(X_t) )\, dt + \alpha(X_t) \, dY_t  + \sigma(X_t) \, dW_t.
\end{align}
We now set
\begin{align}\label{def:rho}
\rho_t f \triangleq \tE \left[ f(X_t) \zeta_t^{-1} \Big| \, \Fy_t \right],
\end{align}
with $\zeta_t$ defined in \eqref{eq:tilde-P}.
Then by Bayes formula, $\pi_t f = \frac{ \rho_t f}{\rho_t 1}$ and moreover, $\rho_t f$ solves
the \emph{linear} stochastic partial differential equation
\begin{align}\label{eq:zakai}
d(\rho_t f) = \rho_t( A f) \, dt + \sum_{k=1}^{d_Y} \left[ \rho_t( h_k f) + \rho_t( B^k f) \right] \,
dY^k_t,
\end{align}
with $A, B^k$ from \eqref{eq:A-B-operators}. The measure-valued Markov process $\rho_t$ is called the \emph{unnormalized} conditional distribution of $X$ and will play a major role in the subsequent analysis. Under the given smoothness assumptions, it is known \cite{BensoussanBook} that $\pi_t$ (and $\rho_t$) will possess a smooth density in $W^1_p$ for all $p > 1$ and $t > 0$.

Returning to our optimal stopping problem \eqref{eq:basic-problem}, let us define the value function $V$ by
$$ V(t, \xi, y; T) \triangleq \sup_{\tau \le T} \E \left[ g(\tau, X_\tau, Y_\tau) \big| \,X_t \sim \xi, Y_t = y  \right].$$
Economically, $V$ denotes the optimal reward that can be obtained on the horizon $[t,T]$ starting with initial condition $Y_t=y$ and $X_t \sim \xi$. Using conditional expectations we may write,
\begin{align}\notag
 V(t, \xi, y) & = \sup_{t \le \tau \le T} \E^{t,\xi,y} \left[ \pi_\tau g(\tau, \cdot, Y_\tau)\right]  \\  \notag 
  & = \sup_{t \le \tau \leq T}  \tE^{t,y,\xi}[ \rho_{\tau} g(\tau, \cdot,Y_\tau)] \\ \label{eq:mainProb} & = \sup_{t \le \tau \le T}  \tE^{t,\xi,y}[ G(\tau, \rho_\tau, Y_\tau)], \qquad\text{where}\quad G(t,\xi,y) \triangleq \int_{\R^d} g(t,x,y) \xi(dx),
\end{align}
and where $\tE^{t,\xi,y}$ denotes $\tP$-expectation conditional on $Y_t=y, \rho_t = \xi$.

Equation \eqref{eq:mainProb} achieved two key transformations. First, its right-hand-side
is now a standard optimal stopping problem featuring the Markov hyperstate $(\rho_t, Y_t)$. Secondly, \eqref{eq:mainProb} has separated the filtering and optimization steps by introducing the fully observed problem through the new state variable $\rho_t$. However, this new formulation remains complex as $\rho_t$  is an infinite-dimensional object. With a slight abuse of notation, we will write $V(t,\rho_t,Y_t)$ to denote the value function as a function of the current unnormalized distribution $\rho_t$.

As can be seen from the last two lines of \eqref{eq:mainProb}, one may solve \eqref{eq:basic-problem} either under the original \emph{physical} measure $\PP$ using $\pi_t$, or equivalently  under the reference measure $\tP$ using $\rho_t$. In our approach we will work with the latter formulation due to the simpler dynamics of $\rho_t$ and more importantly due to the fact that under $\tP$ one can separate the evolution of $Y$ and $X$. In particular, under $\tP$, $Y$ is a Brownian motion and can be simulated entirely on its own. In contrast, under $\PP$, the evolutions of $\pi_t$ and $Y$ are intrinsically tied together due to the \emph{joint} (and unobserved) noise source $(U_t)$.



\subsection{Snell Envelope}
Let us briefly summarize the Snell envelope theory of optimal stopping in our setting. All our results are stated under the $\tP$ reference measure, following the formulation in \eqref{eq:mainProb}.
For any $\Fy$-stopping time $\sigma$, define
\begin{align}
\Z(\sigma) = \sup_{\sigma \leq \tau \leq T} \tE \left[ G(\tau, \rho_\tau, Y_\tau) | \Fy_\sigma \right].
\end{align}

\begin{prop}[\cite{Mazziotto86}]
The set
$(\Z_\sigma)$ form a supermartingale family, i.e.\ there exists a continuous process $\cZ$, such that
$\Z(\sigma) = \cZ_\sigma$, $\cZ$ stopped at time $\sigma$. Moreover, an optimal time $\tau$ for \eqref{eq:mainProb} exists and is given by
$\tau = \inf \{t \colon \cZ_t = G(t, \rho_t, Y_t)\}$.
\end{prop}

The above process $\cZ$ is called the Snell envelope of the optimal stopping problem \eqref{eq:mainProb}. The proposition implies that to solve \eqref{eq:mainProb} it suffices to compute the Snell envelope $\cZ$. We denote by $t \le \tau^*_t \le T$ the optimal $\tau$ achieving the supremum in $\cZ_t = \tE[ G(\tau^*_t, \rho_{\tau^*_t}, Y_{\tau^*_t}) | \Fy_t]$. By virtue of the (strong) Markov property of $(\rho_t, Y_t)$ and the fact that $\rho_t$ is a sufficient statistic for the distribution of $X_t | \Fy_t$ it follows that $V(t, \rho_t,Y_t) = \sup_{t \le \tau \le T} \tE[ G(\tau, \rho_\tau, Y_\tau) | \Fy_t] = \cZ_t$ and \eqref{eq:basic-problem} is equivalent to finding $\tau_0^*$ above. Mazziotto \cite{Mazziotto86} also gave a formal proof of the equivalence of the Snell envelopes under $\PP$ and $\tP$ that we discussed in the end of the previous section.

%

To make computational progress in computing $\tau^*_0$,  it will be eventually necessary to discretize time. Thus, we restrict possible stopping times to lie in the set $\mathcal{S}^\Delta = \{0, \Dt, 2\Dt, \ldots, T\}$, and label the corresponding value function (of the so-called Bermudan problem) as
$$
V^{\Delta}(t,\xi,y) = \sup \{ \tE^{t,\xi,y}[ G(\tau, \rho_\tau, Y_\tau) ] : \tau \text{ is } \mathcal{S}^\Delta\text{-valued  }, \Fy\text{-adapted}\}.$$
 In this discrete version, since one either stops at $t$ or waits till $t+\Dt$, the dynamic programming principle implies that the Snell envelope satisfies
\begin{align}\label{eq:discrete-snell}
V^\Delta(t,\rho_t,Y_t) & = \max \Bigl( G(t,\rho_t, Y_t), \tE[ V^\Delta(t+\Dt,\rho_{t+\Dt},Y_{t+\Dt}) | \Fy_t] \Bigr).
\end{align}

\subsection{Continuation Values and Cashflow Functions}
For notational convenience, we now write $Z_t \equiv (t, \rho_t, Y_t)$ and $G_t = G(Z_t)$. Let $$q_t = q_t(Z_t) \triangleq \tE[ V^\Delta(Z_{t+\Dt}) | \Fy_t],$$ denote the \emph{continuation} value. Then the Snell envelope property \eqref{eq:discrete-snell} implies that
$q_t$ satisfies the recursive equation
\begin{align}\label{eq:cont-value-recursion}
q_t = \tE \left[ \max( G_{t+\Dt}, q_{t+\Dt} ) | \Fy_t \right].
\end{align}
The optimal stopping time $\tau^*_t$ also satisfies a recursion, namely
\begin{align}\label{def:tau}
\tau^*_t = \tau^*_{t+\Dt} 1_{\{q_t > G_t\}} + t 1_{\{q_t \le G_t\}}.
\end{align}
In other words, when the continuation value is bigger than the immediate expected reward, it is optimal to wait; otherwise it is optimal to stop. Equation \eqref{def:tau} also highlights the fact that the continuation value $q_t$ serves as a \emph{threshold} in making the stopping decision. Associated with a stopping rule $\tau^*$ defined above is the future cashflow function.
Denote $B_t(q) \triangleq 1_{\{q_t \le G_t\}}$ and its complement by $B^c_t(q) \equiv 1-B_t(q)$, and starting from the timepoint $t$, define the expected future cashflow as
\begin{align}\label{def:vth}
\vth_t(q)(\bZ) \triangleq \sum_{s=t}^T G(Z_s) 1_{B^c_t \cdot B^c_\tDt \ldots \cdot B_s(q)}.
\end{align}
$\vth_t(q)$ is a path function whose value depends on the realization of $(Z_t)$ between $t$ and $T$, as well as the threshold function $q$. Note that \eqref{def:vth} can be defined for \emph{any} threshold rule $q'$ by simply using $B_t(q')$ instead.
In discrete time using the fact that $\tau^*_t$ is an $\Fy$-stopping time and \eqref{def:tau} 
we get
\begin{align}\notag
q_t(Z_t) = \tE[ V^\Delta(Z_{t+\Dt}) | \Fy_t] & = \tE[ G(\tau^*_{t+\Dt},\rho_{\tau^*_{t+\Dt}}, Y_{\tau^*_{t+\Dt}}) | \Fy_t] \\ \label{eq:q-versus-theta}
& = \tE \left[ \sum_{s=t+\Dt}^T G(s,\rho_s,Y_s) 1_{\{\tau^*_{t+\Dt} = s\}} \Big|\, \Fy_t \right]  =  \tE[ \vth_\tDt(q)(\bZ) | \Fy_t].
\end{align}

It follows that knowing $\vth(q)$, one can back-out the continuation values $q$ and then recover the value function itself from $ V^\Delta(Z_t) = \max( G(Z_t), q(Z_t))$.
In particular, for $t=0$, we obtain $V^\Delta(0,\xi_0,y_0) = \max( G(Z_0), q_0(Z_0))$. The approximation algorithm will compute $q$ and the associated $\vth$ by repeatedly evaluating the conditional expectation in \eqref{eq:q-versus-theta} and updating \eqref{def:vth}. The advantage in using cashflows $\vth(q)$ rather than $q$ itself is that an error in computing $q$ is not propagated backwards \emph{unless} it leads to a wrong stopping decision for \eqref{def:tau}.  As a result, the numerical scheme is more stable.

\begin{rem}
Egloff \cite{Egloff05} discusses a slightly more general situation, where the look-ahead cashflows $\vth$ are taken not on the full horizon $[t,T]$ but only some number $w$ of steps ahead. This then produces
\begin{align}
\vth_{t,w}(q)(\bZ) = \sum_{s=t+\Dt}^{t+w \Dt} G(Z_s) \cdot 1_{B^c_t B^c_\tDt \ldots \cdot B_s} + q_{t+w\Dt}(Z_{t+w\Dt}) \cdot 1_{B^c_t B^c_\tDt \ldots \cdot B^c_{t+w\Dt}},
\end{align}
and one still has $q_t(Z_t) = \tE[ \vth_{t+1,w}(q)(\bZ) | \Fy_t]$ for any $w = 0, \ldots, T-t-1$. In particular, the case $w=0$ is the Tsitsiklis-van Roy \cite{TsitsiklisVanRoy00} algorithm,
\begin{align}\label{eq:tvr}
\vth_{t,0}(q) = G(Z_t) 1_{B_t} + q_t(Z_t)1_{B^c_t}  = \max( G_t, q_t).
\end{align}
\end{rem}

To compute \eqref{eq:q-versus-theta}, the corresponding conditional expectation will be approximated by a finite-dimensional projection $\H$. Indeed, by definition of conditional expectation with respect to the Markov state $(\rho_t, Y_t)$, we have
$q_t(Z_t) = \tE[ \vth_{t+\Dt}(q)(\bZ) | \Fy_t] = F(\rho_t, Y_t)$ for some function $F$. Let $(B_j)_{j=1}^\infty$ be a (Schauder) basis for the Banach space $\R_+ \times \cP(\R^d)$. Then as $r \to \infty$, $F$ (and $q_t)$ can be approximated arbitrarily well by the truncated sum
%
%
%
\begin{align}\label{eq:h-project}
q_t(Z_t) \simeq \hq_t(Z_t) \triangleq \sum_{j=1}^r \alpha_j B_j(\rho_t, Y_t) = \prH \circ \tE[ \vth_{t+\Dt}(q)(\bZ) | \Fy_t],
\end{align}
where the projection manifold (or architecture) is $\H = span( B_j(\xi,y), j=1,\ldots, r)$.  As long
as \eqref{eq:h-project} does not modify much the resulting stopping sets $B_t(\hq)$, one expects that the resulting cashflow function $\vth(\hq)$ will be close to the true one $\vth(q)$. In our filtering context, the extra modification is that $Z_t$ must itself be approximated by a finite-dimensional filter $Z^n_t$. However, if the approximation is high-dimensional, then it should have very little effect on the projection step of the Snell envelope in \eqref{eq:h-project}. 


\subsection{Analytic Approach} The analytic approach to optimal stopping theory characterizes
the value function $V(t,\xi,y)$ in terms of a parabolic-type free boundary problem. This is in
direct counterpart to standard optimal stopping problems for diffusion models.


The major difficulty is the infinite-dimensional nature of the state variable $\pi$. Limited results exist for the corresponding optimal stopping problems on Polish spaces, see e.g.\ \cite{MazziottoStettnerZabczyk88,Mazziotto86}. In particular,  \cite{MazziottoStettnerZabczyk88} characterize $V$ as the minimal excessive function dominating $G$ in terms of the (Feller) transition semigroups of $(\pi_t, Y_t)$. A more direct theory is available when $\pi_t \in H$ belongs to a Hilbert space; this will be the case if $\xi_0$ (and therefore $\pi_t$ for all $t$) admits a smooth $L^2$-density. Even then, since the smoothness properties of $V$ with
respect to $\xi$ are unknown, one must work with viscosity solutions to second-order pdes as
is common in general stochastic control. The following proposition is analogous to Theorem 2.2
in \cite{GatarekSwiech}. Denote by $D$ the Fr\'echet derivative operator and for a twice Fr\'echet differentiable test function $\phi(t,\xi,y)$ let
\begin{align}
\mathcal{L}\phi = \half \tr \left( (\sigma \sigma^T + \alpha \alpha^T) D^2_{\xi\xi} \phi \right) + \langle b, D_\xi\phi \rangle + h \partial_y \phi + \alpha D_{\xi}\phi \cdot \partial_{y}\phi,
\end{align}
with $\langle \cdot, \cdot \rangle$ denoting the inner product in $H$ 
be the infinitesimal generator of the Markov process $(\pi_t,Y_t)$. 

\begin{prop}
The value function $V(t,\pi,y)$ is the unique viscosity solution of
\begin{align}\label{eq:freebnd}
\left\{ \begin{aligned} V_t + \mathcal{L}V \leq 0, \\
V(t,\pi,y) \geq G(t,\pi, y).\end{aligned} \right.
\end{align}
Moreover, $V$ is bounded and locally Lipschitz (with respect to the Hilbert norm).  
\end{prop}
In principle the infinite-dimensional free boundary problem \eqref{eq:freebnd} can be tackled by a variety of numerical methods including the projection approach that passes to a finite-dimensional subset of $L^2(\R^d)$. 



\section{New Algorithm}\label{sec:algorithm}
In this section we describe a new numerical simulation algorithm to solve \eqref{eq:mainProb}. This algorithm will be a combination of the minimal-variance branching particle filter algorithm for
approximating $\pi_t$ and $\rho_t$, described in Section \ref{sec:particle}, and the regression Monte Carlo algorithm described in Section \ref{sec:rms}.

\subsection{Particle Filtering}\label{sec:particle}
The main idea of particle filters is to approximate the measure-valued conditional
distribution $\pi_t$ by a discrete system of point masses that follows a
mutation-selection algorithm to reproduce the dynamics of \eqref{eq:zakai}. In what follows we
summarize the particular algorithm proposed in \cite{Crisan06,CrisanLyons99,Gaines98}. We assume that we are given \eqref{eq:X-sde}-\eqref{eq:Y-sde} with continuous observation of $(Y_t)$.
Fix $n>0$; we shall approximate  $\pi_t$ by a particle system $\pi^n_t$ of $n$ particles. The interacting
particle system consists of a collection of $n$ weights $a^n_j(t)$ and corresponding locations
$v^n_j(t) \in \R^{d}$, $j=1,\ldots, n$. We think of $v^n_j$ as describing the evolution of the $n$-th particle and of $a^n_j(t) \in \R_+$ as its importance in the overall system.  Begin by initializing the system by independently drawing $v^n_j(0)$ from the initial distribution $X_0 \sim \xi_0$ and taking $a^n_j(0) = 1\; \forall j$. Let $\delta$ be a parameter indicating the frequency of mutations; the description below is for a generic time
step $t \in [m \delta, (m+1)\delta)$, assuming that we already have $v^n_j(m \delta)$ and $a^n_j(m \delta) \equiv 1$.

First, for $m \delta \le t < (m+1)\delta$ we have
\begin{align}\label{eq:v-evolution}
v^n_j(t) = v^n_j( m\delta) + \int_{m\delta}^t (b - \alpha h)(v^n_j(s) )\, ds + \int_{m\delta}^t \alpha( v^n_j(s) )\, dY_s + \int_{m\delta}^t \sigma( v^n_j(s) )\, dW^{(j)}_s,
\end{align}
where $W^{(j)}$ are $n$ independent $\tP$-Wiener processes. Thus, each particle location evolves independently according to the law of $X$ under $\tP$. The unnormalized weights $a^n_j(\cdot)$ are given by the stochastic exponentials
\begin{align}\label{eq:a-evolution}
a^n_j(t) = 1 + \sum_{k=1}^{d_Y} \int_{m\delta}^t a^n_j(s) h_k(v^n_j(s))\, dY^k_s = \exp \left( \sum_{k=1}^{d_Y} \int_{m\delta}^t h_k( v^n_j(s) )\, dY^k_s - \half \sum_{k=1}^{d_Y} \int_{m\delta}^t h_k( v^n_j(s))^2 \, ds \right).
\end{align}

Let $$\bar{a}^n_j((m+1)\delta- ) \triangleq \frac{a^n_j( (m+1)\delta-) }{\sum_j a^n_j( (m+1)\delta-)} \in (0,1),$$ denote the normalized weights at the next mutation time. Then at $t=(m+1)\delta$ each particle produces $o^n_j( (m+1)\delta)$ offspring inheriting the parent's location, with the branching carried out such that
\begin{align} \left\{ \begin{aligned}
o^n_j( (m+1)\delta ) & = \left\{ \begin{aligned} \lfloor  n \bar{a}^n_j( (m+1)\delta-) \rfloor & \text{ with prob. } \quad 1- \{ n \bar{a}^n_j( (m+1)\delta- ) \}, \\
1 + \lfloor n \bar{a}^n_j( (m+1)\delta-) \rfloor & \text{ with prob. }\quad \{ n \bar{a}^n_j( (m+1)\delta-) \},
\end{aligned} \right. \\
\sum_{j=1}^n o^n_j( (m+1)\delta ) & = n, \end{aligned} \right.
\end{align}
where $\{ x \}$ denotes the fractional part of $x\in \R$.  Note that the different $o^n_j$'s are correlated so that the total number of particles always stays constant at $n$. One way to generate such $o^n_j$'s is given in the Appendix of \cite{Crisan06}. Following the mutation, particle weights are reset to $a^n_j( (m+1)\delta) = 1$ and one proceeds with the next propagation step.

With this construction we now set for $m\delta \le t < (m+1)\delta$,
\begin{align}\label{eq:particle-filter}
\left\{ \begin{aligned}
\pi^n_t & \defn \sum_{j=1}^n \frac{ n a^n_j(t) }{\sum_{\ell =1}^n a^n_\ell(t)} \delta_{v^n_j(t)}(\cdot);\\
 \rho^n_{t}  & \defn \left[\prod_{\ell=1}^m \left( \frac{1}{n} \sum_{j=1}^n a^n_j(\ell \delta-) \right) \right] \cdot \left( \frac{1}{n} \sum_{j=1}^n a^n_j(t) \delta_{v^n_j(t)} (\cdot)\right).
 \end{aligned} \right.
\end{align}
Interpreted as a probability measure on $\R^{d}$, $\pi^n_t$ ($\rho^n_t)$ is an approximation to the true $\pi_t$ (resp.\ $\rho_t$) as indicated by the following 
\begin{prop}[\cite{Crisan06}, Theorem 5]\label{prop:rho-error}
There exist constants $C_1(t),C_2(t)$ such that for any $f \in \cC^1_b(\R^d)$,
\begin{align}\label{eq:rho-error}
\tE \Bigl[ (\frac{ \rho^n_t f - \rho_t f}{\rho_t 1})^2\Bigr] \leq \frac{C_1(t)}{n} \| f \|^2_{1,\infty},
\end{align}
which in turn implies that (since $\E[ \zeta_t^2]$ is bounded)
\begin{align}\label{eq:pi-error}
\E\Bigl[ (\pi^n_t f - \pi_t f)^2 \Bigr] \le 
\frac{C_2(t)}{n} \| f \|^2_{1,\infty},
\end{align}
with $C_i(t) = \cO(\e^t \cdot t)$.
\end{prop}

Similar results can be obtained under the assumption that $Y$ is observed \emph{discretely}
every $\delta$ time units. In that case one simply takes,
\begin{align*}
a^n_j((m+1)\delta-) = \exp \left( \sum_{k=1}^{d_Y} h_k( v^n_j(m\delta) )\cdot (Y^k_{(m+1)\delta}-Y^k_{m\delta}) - \half \sum_{k=1}^{d_Y} h_k( v^n_j(m\delta))^2 \cdot\delta \right),
\end{align*}
with the rest of the algorithm remaining unchanged.

The use of discrete point masses in the interacting particle filter renders the analytical results based on Hilbert-space theory (e.g.\ \eqref{eq:freebnd}) inapplicable. This can be overcome by considering regularized particle filters \cite{LeGlandOudjane}, where point masses are replaced by smooth continuous distributions and the particle branching procedure switches back to a true re-sampling step.

\subsection{Regression Monte Carlo}\label{sec:rms}
The main idea of our algorithm is to simulate $N$ paths of the $Z$ process (or rather the particle approximation $(Z^n)$), yielding a sample $(z^k_t)$,
$k=1,2,\ldots, N$, $t = 0, \Dt, \ldots, T$. 
To simulate $(z^i_t)$, we first simulate the Brownian motion $(Y_t)$ under $\tP$, and then re-compute $\rho^n_t$ along the simulated paths as described in the previous subsection.
 Using this sample and
approximation architectures $\H_t$ of \eqref{eq:h-project}, we approximate the projection $\prH$ through an empirical least-squares regression. Namely, an empirical continuation value is computed according to
\begin{align}
\hq_t = \argmin_{f \in \H_t} \frac{1}{N} \sum_{i=1}^N | f(z^i_t) - \vth^N(\hq)(z^i_{t+\Dt}) |^2 \simeq \prH \circ \E[ \vth^N_{t+\Dt}(\hq)(\bZ^n) | \Fy_t],
\end{align}
where $\vth^N$ is the empirical cashflow function along simulated paths obtained using the future $\hq$'s. One then updates pathwise $\vth^N$ and $\tau$ using \eqref{def:vth} and \eqref{def:tau} respectively and proceeds recursively backwards in time. This is the same idea as the celebrated regression Monte Carlo algorithm of Longstaff and Schwartz \cite{Longstaff}. The resulting error between $\hq$ and the true $q$ will be studied in Section \ref{sec:errors} below.

Many choices exist regarding the selection of basis functions $B_j(\rho_t, Y_t)$ for the regression step. As a function of $y$, one may pick any basis for $L^2(\R^{d_Y},\tP)$, e.g.\ the Laguerre polynomials. As a function of $\rho$, a natural probabilistic choice involves the moments of $X_t|\Fy_t$, i.e.\ $\sum_i \alpha_i (\rho_t x^i)$. It is also known that using
a basis function of the form $EUR(z) \triangleq \tE_t[ G(Z_T)]$ (the conditional expectation of the terminal reward or the ``European'' counterpart,) is a good empirical choice.

\begin{rem}
If one only uses the first two conditional moments of $X$, $\rho_t x$ and $\rho_t x^2$ inside the basis functions, then our algorithm can be seen as the non-Markovian analogue of applying the Extended Kalman filter for the partial observations of $X$ and then computing the (pseudo)-Snell envelope of \eqref{eq:basic-problem}. In that sense, our approach generalizes all the previous filtering projection methods \cite{BrigoHanzonIEEE98,BudhirajaKushnerIEEE} for \eqref{eq:basic-problem}.
\end{rem}

\subsection{Overall Algorithm}\label{sec:algo-summary}
For the reader's convenience, we now summarize the overall algorithm for solving \eqref{eq:basic-problem}.
\begin{itemize}
\item Select model parameters $N$ (number of paths); $n$ (number of particles per path); $\Dt$ (time step for Snell envelope); $\delta$ (time step for observations and particle mutation); $B_i$ (regression basis functions); $r$ (number of basis functions).

\item Simulate $N$ paths of $(y^k_t)$ under $\tP$ (which is a Brownian motion) with fixed initial condition $y^k_0 = y_0$.

\item  Given the path $(y^k_t)$, use the particle filter algorithm \eqref{eq:v-evolution}-\eqref{eq:a-evolution}-\eqref{eq:particle-filter} to compute $\rho^{n,k}_t$ along that path, starting with $\rho^{n,k}_0 \sim \xi_0$.

\item Initialize $\hq^{k}(T) = \vth^{N,k}_T(\hq) = G(z^k_T)$, $\tau^k(T)=T$, $k=1,\ldots,N$.

\item Repeat for $t=(M-1)\Dt, \ldots, \Dt, 0$:
\begin{itemize}
 \item Evaluate the basis functions $B_i(z^k_t)$, for $i=1,\ldots, r$ and $k=1,\ldots,N$.

 \item Regress $$\alpha^N_t \triangleq \argmin_{\alpha \in \R^{r}}\sum_{k=1}^N  \Bigr|\vth^{N,k}_{t+\Dt}(\hq) - \sum_{i=1}^{r} \alpha^i B_i(  z^k_t) \Bigl|^2.$$

 \item For each $k=1,\ldots,N$ do the following steps: Set $ \hq^k(t) = \sum_{i=1}^{r} \alpha^{N,i}_t B_i( z^k_t)$.

 \item Compute $G(z^k_t) = \rho^{n,k}_t g(t, \cdot, y^k_t)$.

 \item Update $ \vth^{N,k}_t(\hq) = \left\{ \begin{aligned} G(z^k_t) & \quad \text{if }
 \hq^k_t < G(z^k_t); \\
  \vth^{N,k}_{t+\Dt}(\hq)  & \quad\text{otherwise.} \end{aligned}\right.$ 

  \item Update $\tau^k(t) = \left\{ \begin{aligned} t \quad & \quad\text{if } \hq^k(t) < G(z^k_t); \\
       \tau^k(t+\Dt) & \quad\text{otherwise.} \end{aligned}\right.$

\end{itemize}
  \item End loop;

  \item Return $V^\Delta(0, \xi_0, y_0) \simeq \frac{1}{N} \sum_{k=1}^N \vth^{N,k}_0(\hq)$.
  \end{itemize}

Note that it is not necessary to save the entire particle systems $(v^{n,k}_j(m\Dt))_{j=1}^n$ after the simulation step; rather one needs to keep around just the evaluated basis functions $ (B_i(z^k_t))_{i=1}^{r}$, so that the total memory requirements are $\cO( N \cdot M \cdot r)$. In terms of number of operations the overall algorithm complexity  is $\cO( M \cdot N \cdot (n^2 + {r}^3) )$, with the most intensive steps being the resampling of the filter particles and the regression step against the $r$ basis functions.






\section{Error Analysis}\label{sec:errors}
This section is devoted to the error analysis of the algorithm proposed in Section \ref{sec:algo-summary}. Looking back, our numerical scheme involves three main errors. These are:
\begin{itemize}
\item Error in computing $\rho_t$ which arises from using a finite number of particles and the resampling error of the particle filter $\rho^n_t$;

\item Error in projecting the cashflow function $\vth$ onto the span of basis functions $\H$ and the subsequent wrong stopping decisions;

\item Error in computing projection coefficients $\alpha^i$ due to the use of finite-sample least-squares regression.
\end{itemize}

We note that the filtering error is propagated forward, while the projection and empirical errors are propagated backwards. In that sense, the filtering error is more severe and should be controlled tightly. The projection error is the most difficult to deal with since we only have crude estimates on the dependence of the value function on $\rho_t$. Consequently, the provable error estimates are very pessimistic. Heuristic considerations would imply that this error is in fact likely to be small. Indeed, the approximate decision rule will be excellent as long as $\tP( \{ q_t > G_t \} \cap \{ \hq_t \le G_t \} )$ is small, since the given event is the only way that the optimal cashflows are computed incorrectly. By applying domain knowledge the above probability can be controlled through customizing the projection architecture $\H_t$. For instance, as mentioned above, using $EUR(z)$ as one of the basis functions is often useful.

The sample regression error is compounded due to the fact that we do not use the true basis functions but rather approximations based on $Z^n$. This implies the presence of \emph{error-in-variable} during the regression step from the pathwise filtering errors.  It is well-known (see e.g.\ \cite{Fuller87}) that this leads to attenuation in the computed regression result, i.e.\ $|\alpha^{N,i}| \le |\alpha^i|$. An extensive statistical literature treats error reduction methods to counteract this effect, a topic that we leave to future research.


%
%
%


As a notational shorthand, in the remainder of this section we write $\tE_t$ to denote expectations (as a function on $\R^{d} \times \R^{d_Y}$) conditional on $Y_t=y$ and $\rho_t = \xi$. We recall that the optimal cashflows satisfy
$$ q_t = \tE_t [\vth_\tDt(q)(\bZ)], $$ while the approximate cashflows are
$$ \hq_t =  \prH^N \circ \tE_t[\vth^N_\tDt(\hq)(\bZ^n)].$$
Note that inside the algorithm, $\hq_t$ is evaluated not at the true value $Z_t = (\rho_t, Y_t)$, but at the approximate point $Z^n_t$. To emphasize the process under consideration we denote by $q^n_t \equiv q^n_t(Z^n_t)$ the continuation function resulting from working with the $Z^n$-process. Observe that the difference between $q^n$ and the true $q$ is solely
due to the inaccurate recursive evaluation of the reward $G$ (since $Y$ is simulated exactly); thus if the original reward $g$ in \eqref{eq:basic-problem} is independent of $X$ then $q^n \equiv q$.

The error analysis will be undertaken in two steps. In the first step, we consider the mean-squared error between the continuation value $q_t$ based on the true filter $\rho_t$ and the continuation value $q^n_t$ based on the approximate filter $\rho^n_t$. In the second step, we will study the difference between  $q^n_t$ and the approximate $\hq_t$ above. Throughout this section, $\| \cdot \|_2 \equiv \tE[ |\cdot|^2 ]^{1/2}$. 


\begin{lemma}\label{lem:scE-5}
There exists a constant $C(T)$, such that for all $t \le T$
\begin{align}\label{eq:scE-5}
\left\| \tE_t[ \vth_{t+\Dt}(q^n) - \vth_{t+\Dt}(q)] \right\|_2 & \le \frac{(T-t)\cdot C(T)}{\Dt \cdot \sqrt{n}} \cdot \| g \|_{1,\infty}.
\end{align}
\end{lemma}

\begin{proof}
Suppose without loss of generality that $q^n(Z^n_t) > q(Z_t)$. 
Let $\tau$ be an optimal stopping time for the problem represented by $q^n$. Clearly such $\tau$ is sub-optimal for $q$; moreover since both $Z$ and $Z^n$ are $\Fy$-adapted, $\tau$ is admissible for $q$. Therefore,
\begin{align*}
(q^n(Z^n_t) - q(Z_t) )^2 & \le \tE_t \bigl[ G(Z^n_\tau) - G(Z_\tau) \bigr]^2 \\
& = \left\{ \sum_{s=t+\Dt}^T \tE_t \left[ ( G(Z^n_s)-G(Z_s)) \cdot 1_{\{\tau = s\}} \right] \right\}^2 \\
& \le  \sum_{s=t+\Dt}^T \frac{T-t}{\Dt} \cdot \tE_t \bigl[ | G(Z^n_s)-G(Z_s)|^2 \bigr],
\end{align*}
where the last line is due to Jensen's inequality. Averaging over the realizations of $(Z^n_t, Z_t)$ we then obtain
\begin{align*}
\tE[ |q^n(Z^n_t) - q(Z_t) |^2] & \le \sum_{s=t+\Dt}^T \frac{T-t}{\Dt} \cdot \tE[ | G(Z^n_s)-G(Z_s)|^2 ] \\
& \le  \sum_{s=t+\Dt}^T \frac{(T-t)C(T)}{\Dt \cdot n} \|g \|^2_{1,\infty}  = \frac{(T-t)^2 \cdot C(T)}{\Dt^2 \cdot n} \| g \|^2_{1,\infty},
\end{align*}
using Proposition \ref{prop:rho-error}.

Note that this error explodes as $\Dt \to 0$ due to the fact that we do not have tight bounds for $\tE_t[ | G(Z^n_s)-G(Z_s)|^2 1_{\{\tau = s\}}]$. In general, one expects that $\tE_t[ | G(Z^n_s)-G(Z_s)|^2 1_{\{\tau = s\}}] \simeq \tE_t[ | G(Z^n_s)-G(Z_s)|^2 ] \cdot \PP( \tau=s)$ which would eliminate the $\Dt^{-2}$ term on the last line above.

\end{proof}

In the second step we study the $L^2$-difference of the unnormalized continuation values, $\| q^n_t - \hq_t \|_2 \equiv \tE[ (q^n_t(Z^n_t) - \hq_t(Z_t))^2]^{1/2}$. This total error can be decomposed as
\begin{align}\notag
\| \hq_t - q_t \|_2 & \le 
\underbrace{\bigl\| \prH^N \circ \tE_t[ \vth^N_\tDt(\hq)(\bZ^n)] - \prH \circ \tE_t[\vth_\tDt(\hq)(\bZ^n) \bigr\|_2}_{\scE_1}  \\
\label{eq:error-decomp} & \quad + \underbrace{\bigl\| \prH \circ \tE_t[ \vth_\tDt(\hq)(\bZ^n)] - \tE_t [ \vth_\tDt(\hq)(\bZ^n)] \bigr\|_2}_{\scE_2}  + \underbrace{\bigl\| \tE_t[ \vth_\tDt(\hq)(\bZ^n) - \vth_\tDt(q)(\bZ^n)] \bigr\|_2}_{\scE_3}.  
\end{align}
The three error terms $\scE_i$ on the right-hand-side of \eqref{eq:error-decomp} are
respectively 
the empirical error $\scE_1$, the projection
error $\scE_2$, and the recursive error from the next time step $\scE_3$.  
Each of these terms is considered in turn in the next several lemmas with the final result summarized in Theorem \ref{thm:error-decomp}.  The first two lemmas have essentially appeared in
\cite{Egloff05} and the proofs below are provided for completeness.


\begin{lemma}[{\cite[Lemma 6.3]{Egloff05}}]\label{lem:scE-2}
Define the centered loss random variable
\begin{align}\label{eq:centered-loss}
\ell_t(\hq)(\bZ^n) = | \hq_t - \vth_\tDt(\hq)|^2 - | \prH \circ \tE_t [ \vth_\tDt(\hq)] - \vth_\tDt(\hq)|^2.
\end{align}
Then
\begin{align}\label{eq:Lemma-Loss}
\scE_1^2 = \| \hq_t - \prH \circ \tE_t[ \vth_\tDt(\hq)] \|_2^2 \le \tE[ \ell_t(\hq)(\bZ^n) ].
\end{align}
\end{lemma}
\begin{proof}
First note that
\begin{align}\label{eq:egloff-1}
\| \hq_t -  \prH \circ \tE_t[ \vth_\tDt(\hq)] \|^2_2 + \|  \prH \circ \tE_t[ \vth_\tDt(\hq)] -   \tE_t[ \vth_\tDt(\hq)] \|^2_2 \le \| \hq_t - \tE_t[ \vth_\tDt(\hq) ] \|^2_2,
\end{align}
because $\hq_t \in \H_t$ belongs to the convex space $\H_t$, while
$\prH \circ \tE_t[ \vth_\tDt(\hq)] \in \H_t$ is the projection of $\vth_\tDt(\hq)$. Therefore the three respective vectors form an obtuse triangle in $L^2$:
\begin{align*}
\tE\left[ (\hq_t -  \prH \circ \tE_t[ \vth_\tDt(\hq)]) \cdot (\prH \circ \tE_t[ \vth_\tDt(\hq)] -   \tE_t[ \vth_\tDt(\hq)] )\right] \le 0.
\end{align*}

Direct expansion using the tower property of conditional expectations and the fact that $\hq_t \in \Fy_t$ shows that $\tE [ (\hq_t - \tE_t[  \vth_\tDt(\hq)] )\cdot( \tE_t[  \vth_\tDt(\hq)] - \vth_\tDt(\hq) )] = 0,$ so that
\begin{align}\label{eq:egloff-2}
\tE \left[ |\hq_t - \tE_t[  \vth_\tDt(\hq)] |^2 \right] + \tE \left[ | \tE_t[  \vth_\tDt(\hq)] - \vth_\tDt(\hq) |^2 \right] = \tE \left[ | \hq_t - \vth_\tDt(\hq) |^2 \right].
\end{align}
Similarly,
$$ \tE\left[ (\prH \circ \tE_t[ \vth_\tDt(\hq)] -  \tE_t[ \vth_\tDt(\hq)]) \cdot(
\tE_t[ \vth_\tDt(\hq)] - \vth_\tDt(\hq) )\right] = 0,$$
and so
\begin{multline}\label{eq:egloff-3}
\tE \left[|  \prH \circ \tE_t[ \vth_\tDt(\hq)] -  \tE_t[ \vth_\tDt(\hq)] |^2 \right] = \tE \left[|  \prH \circ \tE_t[ \vth_\tDt(\hq)] -  \vth_\tDt(\hq) |^2 \right] \\ - \tE \left[|  \tE_t[ \vth_\tDt(\hq)] -  \vth_\tDt(\hq) |^2 \right].
\end{multline}
Combining \eqref{eq:egloff-1}-\eqref{eq:egloff-2}-\eqref{eq:egloff-3} we find
\begin{align*}
\Bigl\| \hq_t - \prH & \circ \tE_t[ \vth_\tDt(\hq)] \Bigr\|_2^2  \le \tE \Bigl[  |\hq_t -  \vth_\tDt(\hq)|^2 - | \tE_t[  \vth_\tDt(\hq)] -  \vth_\tDt(\hq)|^2 \\  - & \left\{|\prH \circ \tE_t[ \vth_\tDt(\hq)] -  \vth_\tDt(\hq)|^2 -  |\tE_t[ \vth_\tDt(\hq)] -  \vth_\tDt(\hq)|^2\right\} \Bigr] = \tE[ \ell_t(\hq)(\bZ^n)].
\end{align*}
\end{proof}
The above lemma shows that the squared error $\scE_1^2$ resulting from the empirical regression  used to obtain $\hq_t$ (which recall is a proxy for $\tE_t [ \vth_\tDt(\hq)]$) can be expressed as the difference between the expected actual difference $| \hq_t - \vth_\tDt(\hq)|^2$ versus the theoretical best average error after the projection $| \prH \circ \tE_t [ \vth_\tDt(\hq)] - \vth_\tDt(\hq)|^2$.

%
%
%
%
%

\begin{lemma}[cf.\ {\cite[Proposition 6.1]{Egloff05}}]\label{lem:scE-3}
We have $\scE_2 \le 2 \| \tE_t[ \vth_\tDt(q^n) -  \vth_\tDt(\hq)] \|_2 + \inf_{f \in \H_t} \| f - q^n_t \|_2$.
\end{lemma}
\begin{proof}
We re-write,
\begin{align*}
\scE_2 = \left\| \prH \circ \tE_t[ \vth_\tDt(\hq)] - \tE_t [ \vth_\tDt(\hq)] \right\|_2 & \le \left\| \prH \circ \tE_t[ \vth_\tDt(\hq)] - \prH \circ \tE_t [ \vth_\tDt(q^n)] \right\|_2  \\ +  \bigl\| \prH \circ \tE_t [ \vth_\tDt(q^n)] & - \tE_t[ \vth_\tDt(q^n)]  \bigr\|_2 + \| \tE_t[ \vth_\tDt(q^n) -  \vth_\tDt(\hq)] \|_2 \\
& \le 2 \| \tE_t[ \vth_\tDt(q^n) -  \vth_\tDt(\hq)] \|_2 + \inf_{f \in \H_t} \| f - \tE_t[ \vth_\tDt(q^n)]\|_2 \\
& = 2 \left\| \tE_t[ \vth_\tDt(q^n) -  \vth_\tDt(\hq)] \right\|_2 + \inf_{f \in \H_t} \| f - q^n_t \|_2,
\end{align*}
where the second inequality uses the contraction property of the projection map $\prH$ and the definition of projection onto the manifold $\H_t$.
\end{proof}

%

\begin{lemma}\label{lem:scE-4}
We have for any $p > 1$
\begin{align}\label{eq:Lemma-CLP}
\Bigl\| \tE_t[ \vth_{t+\Dt}(q^n) - \vth_{t+\Dt}(\hq)] \Bigr\|_p \le \sum_{s={t+\Dt}}^T \bigl\| \hq_s - q^n_s \bigr\|_p. 
\end{align}
\end{lemma}
\begin{proof}
To simplify notation we drop the function arguments and also write $q_{t+1}, G_{t+1}$, etc., to mean $q_{t+\Dt}$, etc.\ in the proof below.
By definition of the cashflow function, $\scE_3 := \| \tE_t[ \vth_{t+1}(q^n) - \vth_{t+1}(\hq)] \|_p=$
\begin{align*}
 & \quad \Bigl\| \tE_t[ G_{t+1}  1_{
\{q^n_{t+1}  \le G_{t+1} \}} + \vth_{t+2}(q^n)1_{\{q^n_{t+1}  > G_{t+1} \}} - G_{t+1}  1_{\hq_{t+1}  \le G_{t+1} } -\vth_{t+2}(\hq) 1_{\{\hq_{t+1}  > G_{t+1}  \}}]\Bigr\|_p\\
& \quad = \Bigl\| \tE_t[ G_{t+1}  (1_{\{G_{t+1}  \ge q^n_{t+1} \}}-1_{\{G_{t+1}  \ge \hq_{t+1} \}})
+ \vth_{t+2}(q^n)1_{\{q^n_{t+1}  > G_{t+1} \}} -\vth_{t+2}(\hq) 1_{\{\hq_{t+1}  > G_{t+1}  \}}] \Bigr\|_p 
\\
& \quad\le \| \tE_t[A_1] \|_p  +  \Bigl\| \tE_t[ q^n_{t+1} (1_{\{G_{t+1}  \ge q^n_{t+1} \}}-1_{\{G_{t+1}  \ge \hq_{t+1} \}}) + \vth_{t+2}(q^n) 1_{\{G_{t+1}  < q^n_{t+1} \}} - \vth_{t+2}(\hq)1_{\{G_{t+1}  < \hq_{t+1} \}}] \Bigr\|_p,
\end{align*}
where
\begin{align*}
A_1 &= (G_{t+1}  -q^n_{t+1} ) \cdot\left(1_{\{G_{t+1}  \ge q^n_{t+1} \}}-1_{\{G_{t+1}  \ge \hq_{t+1} \}} \right) \\
 & = (G_{t+1}  -q^n_{t+1} ) \left(1_{\{\hq_{t+1}  > G_{t+1}  \ge q^n_{t+1} \}} - 1_{\{ q^n_{t+1}  > G_{t+1}  \ge \hq_{t+1} \}} \right) \\
 & \le (\hq_{t+1} - q^n_{t+1} )1_{\{\hq_{t+1}  > G_{t+1}  \ge q^n_{t+1} \}} - (\hq_{t+1} -q^n_{t+1} )1_{\{ q^n_{t+1}  > G_{t+1}  \ge \hq_{t+1} \}} \\
 & \le | \hq_{t+1} - q^n_{t+1}  |.
\end{align*}
For the remaining terms, using the fact that $q^n_{t+1}  = \tE [ \vth_{t+2}(q^n) | \Fy_{t+1}]$ we obtain
\begin{align*}
\tE_t \left[q^n_{t+1} (1_{\{G_{t+1}  \ge q^n_{t+1} \}}-1_{\{G_{t+1}  \ge \hq_{t+1} \}}) \right] & = \tE_t \left[ \vth_{t+2}(q^n)(1_{\{G_{t+1}  \ge q^n_{t+1} \}}-1_{\{G_{t+1}  \ge \hq_{t+1} \}}) \right],
\end{align*}
and therefore
\begin{align*}
\left\| \tE_t[ \vth_{t+1}(q^n) - \vth_{t+1}(\hq)] \right\|_p & \le \Bigl\| \tE_t \Bigl[ \vth_{t+2}(q^n) \left(1_{\{G_{t+1}  < q^n_{t+1} \}}+1_{\{G_{t+1}  \ge q^n_{t+1} \}}-1_{\{G_{t+1}  \ge \hq_{t+1} \}}\right)  \\
& \qquad\qquad - \vth_{t+2}(\hq)1_{\{G_{t+1}  < \hq_{t+1} \}} \Bigr] \Bigr\|_p + \left\| \tE_t [ | \hq_{t+1} - q^n_{t+1}  | ] \right\|_p\\
& \le \| \hq_{t+1}  - q^n_{t+1}  \|_p + \bigl\| \tE_t[ (\vth_{t+2}(q^n)-\vth_{t+2}(\hq))1_{\{G_{t+1}  < \hq_{t+1} \}} ] \bigr\|_p \\
& \le \| \hq_{t+1}  - q^n_{t+1}  \|_p +\left\| \vth_{t+2}(q^n)-\vth_{t+2}(\hq) \right\|_p.
\end{align*}
By induction, $\scE_3 \le \sum_{s={t+1}}^T \| \hq_s - q^n_s \|_p$ follows.
\end{proof}

Based on Lemmas \ref{lem:scE-5}-\ref{lem:scE-2}-\ref{lem:scE-3}-\ref{lem:scE-4}, we obtain the main
\begin{thm}\label{thm:error-decomp}
We have
\begin{align}\label{eq:thm-1} \| \hq_t(Z^n_t) - q_t(Z_t)\|_2 \le 4^{(T-t)/\Dt} \max_{t
\le s \le T} \left\{ \inf_{f \in \H_s} \| f -  q^n_s \|_2 + \sqrt{\tE[ l_s(\hq)]} \right\} + \frac{C(T)(T-t)}{\Dt \cdot \sqrt{n}} \| g \|_{1,\infty}.\end{align}
\end{thm}

\begin{proof}
Combining Lemmas \ref{lem:scE-2}-\ref{lem:scE-3}-\ref{lem:scE-4} we find that
\begin{align*}
\| \hq_t - q^n_t \|_2 \le \sqrt{\tE[ l_t(\hq)]} + \inf_{f \in \H_t} \| f - q^n_t \|_2 + 3 \sum_{s=t+\Dt}^T \| \hq_s - q^n_s \|_2.
\end{align*}
Therefore, iterating
\begin{align*}
\| \hq_t - q^n_t \|_2 & \le \sqrt{\tE[ l_t(\hq)]} + \inf_{f \in \H_t} \| f - q^n_t \|_2 + 3 \cdot \left(\sqrt{\tE[ l_{t+\Dt}(\hq)]} + \inf_{f \in \H_{t+\Dt}} \| f - q^n_{t+\Dt} \|_2 \right) \\
& \qquad + 9 \cdot \left(\sqrt{\tE[ l_{t+2\Dt}(\hq)]} + \inf_{f \in \H_{t+2\Dt}} \| f - q^n_{t+2\Dt} \|_2 \right) + \ldots \\
& \le 4^{(T-t)/\Dt} \max_{t \le s \le T} \left\{ \sqrt{\tE[ l_s(\hq)]} + \inf_{f \in \H_s} \| f -  q^n_s \|_2 \right\}. 
\end{align*}
Finally, we have 
$ \| \hq_t -q_t \|_2 \le \| \hq_t - q^n_t \|_2 + \| q^n_t - q_t \|_2,$ and applying Lemma \ref{lem:scE-5} the result \eqref{eq:thm-1} follows.
\end{proof}

\subsection{Convergence}
To obtain convergence, one proceeds as follows. First, taking $n\to\infty$ eliminates the filtering error so that $\bZ^n \to \bZ$ and the corresponding errors in evaluating $G$ vanish. Next, one takes $N \to \infty$, reducing the empirical error and the respective centered loss term $\tE[ l_t(\hq)]$. Thirdly, one increases the number of basis functions $r\to\infty$ in order to eliminate the projection error $\inf_{f \in \H_s}\| f - q^n_s \|$.  Finally, taking $\Dt \to 0$ we remove the Snell envelope discretization error.

%

The performed error analysis shows the major trade-off regarding the approximation architectures $\H_t$. On the one hand, $\H_t$ should be large in order to minimize the projection errors $\min_{f \in \H_t} \| f - q^n_t \|$. On the other hand, $\H_t$ should be small to control the empirical variance of the regression coefficients. With many basis functions, one requires a very large number of paths to ensure that $\hq$ is close to $q$. Finally, $\H_t$ should be smooth in order to further bound the empirical regression errors and the filtering error-in-variable accumulated when computing the regression coefficients.

In the original finite-dimensional study of \cite{Egloff05}, the size of $\H_t$ was described in terms of the Vapnik-Cervonenkis (VC) dimensions $n_{VC}$ and the corresponding covering numbers. Using this theory, \cite{Egloff05} showed that overall convergence can be obtained for example by using the polynomial basis for $\H_t$ and taking the number of paths as $N=r^{d+2k}$ where $r$ is the number of basis functions, $d$ is the dimension of the state variable and $k$ is the smoothness of the payoff function $g \in W^k_p$.
In the infinite-dimensional setting of our model, the VC-dimension is meaningless and therefore such estimates do not apply. One could trivially treat $\rho^n_t$ as an $n$-dimensional object, but then the resulting bounds are absurdly poor. It appears difficult to state a useful result on the required relationship between the number of basis functions and the number of paths needed for convergence.

\begin{rem}
A possible alternative is to apply the Tsitsiklis-van Roy algorithm \cite{TsitsiklisVanRoy00}, which directly approximates $q_t$ (rather than $\vth$) using the recursion formula \eqref{eq:tvr}:
$
q_t = \tE_t[ \max(G(Z_\tDt), q(Z_\tDt))].$
 Like in Section \ref{sec:algorithm}, the approximate algorithm would consist in computing via regression Monte Carlo the empirical continuation value
$$
\hq_t = \prH^N \circ \tE_t[ \max(G(Z^n_\tDt), \hq(Z^n_\tDt))].$$
In such a case, the error between $\hq_t$ and $q_t$ admits the simpler decomposition (using $\max(a,b) \le a+b$)
\begin{align}\label{eq:error-tvr}
\left\| \hq_t(Z^n_t) - q_t(Z_t) \right\|_2 & \le \left\| \prH^N \circ \tE_t[ \max(G(Z^n_\tDt), \hq(Z^n_\tDt)] - \prH \circ \tE_t[ \max(G(Z^n_\tDt), \hq(Z^n_\tDt))] \right\|_2 \\ \notag
& \qquad + \left\| \prH \circ \tE_t[ \max(G(Z^n_\tDt), \hq(Z^n_\tDt))] - \tE_t[ \max(G(Z^n_\tDt), \hq(Z^n_\tDt))] \right\|_2 \\ \notag & \qquad +
\left\|  \tE_t[G(Z^n_\tDt) - G(Z_\tDt)] \right\|_2  + \left\| \tE_t[\hq(Z^n_\tDt) - q(Z^n_\tDt)] \right\|_2 \\ \notag & \qquad + \left\| \tE_t[q(Z^n_\tDt) - q(Z_\tDt)] \right\|_2.
\end{align}
We identify the first two terms as the empirical $\scE_1$ and projection $\scE_2$ errors (as in Lemmas \ref{lem:scE-2} and \ref{lem:scE-3}), the third term as the $G$-evaluation error, the fourth term as the next-step recursive error, and finally the last term as the sensitivity error of $q$ with respect to $Z$. Controlling the latter error requires understanding the properties of the continuation (or value) function in terms of current state. This seems difficult in our infinite-dimensional setting and is left to future work. Nevertheless, 
%
%
proceeding as in the previous subsection and iterating \eqref{eq:error-tvr}, we obtain for some constants $C_3,C_4$
$$
\| \hq_t - q_t \|_2 \le \frac{C_3 \cdot (T-t)}{\Dt} \cdot \max_{t \le s \le T} \left\{ \inf_{f \in \H_s} \| f - \hq_s \|_2 + \sqrt{\tE[ l_s(\hq)]} + \frac{C_4}{\sqrt{n}}\| g\|_{1,\infty} + \| q(Z^n_s) - q(Z_s) \|_2 \right\},$$ 
so that the total error is \emph{linear} rather than exponential in number of steps $T/\Dt$ as in Theorem \ref{thm:error-decomp}. Even though this theoretical result appears to be better, empirical evidence shows that the original algorithm is more stable thanks to its use of $\vth$.
\end{rem}

%


\section{Examples}\label{sec:examples}
To illustrate the ideas of Section \ref{sec:algorithm} and to benchmark the described algorithm, 
we consider a model where an explicit finite-dimensional solution is possible.
Let
\begin{align}\label{eq:model-1}
\left\{ \begin{aligned}
dX_t & = -\kappa X_t \, dt + \sigma_X (\rho \, dW_t + \sqrt{1-\rho^2} \, dU_t); \\
dY_t & = (X_t - a)\, dt + \sigma_Y \, dW_t,  \end{aligned} \right. 
\end{align}
with $(U, W)$ being two standard independent one-dimensional Brownian motions.
Thus, $Y$ is a linear diffusion with a stochastic, zero-mean-reverting Gaussian drift $X$. We study the finite horizon optimal stopping problem of the form
\begin{align}\label{eq:obj-1}
V(t,\xi,y) = \sup_{\tau \le T} \E^{t,\xi,y} [\e^{-r \tau}g(X_\tau, Y_\tau)] \triangleq \sup_{\tau \le T} \E^{t,\xi,y} \left[ \e^{-r \tau} (Y_\tau(c_1+X_\tau)-c_2)_+ \right], \qquad c_i \in \R,
\end{align}
which can be viewed as an exotic Call option on $Y$, see the first example in Section \ref{sec:apps}. Note that the payoff is guaranteed to be non-negative even if the controller stops when $Y_\tau(c_1 + X_\tau) < c_2$. In this example, under the reference measure $\tP$, we have
\begin{align*} \left\{ \begin{aligned}
dY_t & = \sigma_Y d\bW_t; \\
dX_t & = [-\kappa X_t - \rho (\sigma_X/\sigma_Y) \cdot (X_t-a) ] \, dt + \rho \sigma_X \,d\bW_t + \sqrt{1-\rho^2} \sigma_X \, dW^\perp_t; \\
d\rho_t(x) & = \frac{1}{2} \sigma_X^2 \rho_t''(x) + \kappa x \rho_t'(x) + \kappa \rho_t(x) + [\frac{x-a}{\sigma_Y^2} - \frac{\rho \sigma_X}{\sigma_Y}]\,dY_t,\end{aligned} \right.
\end{align*}
where $W^\perp$ is a $\tP$-Wiener process independent of $\bW$.

Below we carry out a numerical study with parameter values taken as $$ \begin{array}{c|ccccccccc}
\text{Parameter} & \kappa & a  & \sigma_Y & \sigma_X & T & r & \rho & c_1 & c_2\\ \hline
\text{Value} & 2 & 0.05 & 0.1 & 0.3 & 1 & 0.1 & 0.6 & 1 & 2\\
\end{array}.$$
Since on average $X_t$ is around $0 < a$, $Y$ tends to decrease, so that in \eqref{eq:obj-1} it is optimal to stop early. However, the drift process $X$ is highly volatile and quite often $X_t
> a$ produces positive drift for $Y$, in which case one should wait. Consequently, the stopping region will be highly sensitive to the conditional distribution $\pi_t$.

\subsection{Kalman Filter Formulation}
The model \eqref{eq:model-1} also fits into the Kalman-Bucy \cite{KalmanBucy} filter framework. Thus, if the initial
distribution $X_0 \sim \mathcal{N}(m_0, P_0)$ is a Gaussian density, then $X_t | \Fy_t\sim \mathcal{N}(m_t, P_t)$ is conditionally Gaussian, where \begin{align}\label{eq:kalman-filter-equations}
\left\{
\begin{aligned}
dm_t &= -\kappa m_t\, dt + (\rho \sigma_X + P_t/\sigma_Y) \, d\bW_t, \qquad\qquad\qquad d\bW_t = \frac{dY_t - (m_t-a)\,dt}{\sigma_Y}, \\
dP_t &= ( -2\kappa P_t + \sigma_X^2 - \left( \rho \sigma_X + P_t/\sigma_Y \right)^2 )\,dt.
\end{aligned}\right.
\end{align}
Note that the conditional variance $P_t$ is deterministic and solves the Riccati ode on the second line of \eqref{eq:kalman-filter-equations}.
In \eqref{eq:kalman-filter-equations}, $\bW$ is a
$\PP$-Brownian motion, the so-called innovation process. Moreover, as shown by \cite[Section
12.1]{LiptserShiryaevII}, $\Fy_t = \mathcal{F}^{\bW,\xi_0}_t$, so that we may equivalently
write
\begin{align*}
dY_t = (m_t - a)\,dt + \sigma_Y d\bW_t.
\end{align*}
The pair $(m_t, P_t)$ are sufficient statistics for the conditional distribution of $X_t | \Fy_t$ and the corresponding payoff can be computed as
\begin{align*}
\E[ g(X_t,Y_t) | \Fy_t ] 
& = \E\left[ (y(c_1+m_t+\sqrt{P_t}\mathcal{X})-c_2)_+\right], \qquad\qquad \text{where}\quad \mathcal{X} \sim \mathcal{N}(0,1) \\
& 
= \int_{x^*}^\infty \frac{1}{\sqrt{2 \pi}} \e^{-x^2/2} \left\{( (c_1+m_t) y -c_2) + y \sqrt{P_t} x \right\} dx \\
& = \frac{y \sqrt{P_t}}{\sqrt{2\pi}}\cdot \e^{-(x^*)^2/2} + ( (c_1 + m_t) y -c_2) \cdot(1-\Phi(x^*)) =: G(m_t,P_t,Y_t),
\end{align*}
where $x^* = \frac{c_2 - (c_1+ m_t) y}{y \sqrt{P_t}}$, and $\Phi(x)$ is the standard normal cumulative distribution function. Thus, the original problem is reduced to 
\begin{align}\label{eq:ex1-v}
V(t,m,p,y) = \sup_{\tau \le T} \E \left[ \e^{-r\tau} G(m_\tau, P_\tau, Y_\tau) \Big| \, m_0 = m, P_0=p, Y_0 = y \right].
\end{align}
This two-dimensional problem (recall that $(P_t)$ is deterministic) can be solved numerically using a pde solver applied to the corresponding version of the free boundary problem \eqref{eq:freebnd}. Namely $V$ of \eqref{eq:ex1-v} is characterized by the quasi-variational inequality
\begin{align}\label{eq:example1-pde}
\left\{ \begin{aligned} & \max \Bigl\{
V_t  + (m-a)V_y + \half \sigma_Y^2 V_{yy} - \kappa m V_m + \half (\rho \sigma_X +
P_t/\sigma_Y)^2 V_{mm} \\ & \qquad\qquad + (\rho \sigma_X \sigma_Y + P_t) V_{my} - rV, \, G(m,p,y) - V(t,m,p,y) \Bigr\} = 0,\\ &
 V(T,m,p,y) = G(m,p,y). \end{aligned}\right.
\end{align}

\subsection{Numerical Results}
To benchmark the proposed algorithm we proceed to compare two solutions of \eqref{eq:obj-1}, namely (i) a simulation algorithm of Section \ref{sec:algo-summary} and (ii) a finite-differences pde solver of \eqref{eq:example1-pde}. The Monte Carlo implementation used $N=30\,000$ paths with $n=500$, $\delta = 0.01$, $\Dt = 0.05$ or twenty time-steps.  For basis functions we used the set $\{ 1, y, y^2, \rho_t x, \rho_t g, \rho_t EUR\}$, where $EUR(t,\xi,y) = \tE^{t,\xi,y}[ \e^{-r(T-t)}g(X_T, Y_T)]$ is the conditional expectation of terminal payoff. A straightforward code written in Matlab with minimal optimization took about three minutes to run on a desktop PC. The pde solver utilized a basic explicit scheme and used a $400 \times 400$ grid with $8000$ timesteps. In order to allow a fair comparison, the pde solver also allowed only $T/\Dt =20$ exercise opportunities by enforcing the barrier condition $V(t,m,p,y) \ge G(m,p,y)$ only for $t =m\Dt$, $ m=0,1,\ldots 20$. In financial lingo, we thus studied the Bermudan variant of \eqref{eq:obj-1} with $\Dt = 0.05$.

The obtained results are summarized in Table \ref{table:ex1} for a variety of initial conditions $(\xi_0, Y_0)$. Using the pde solver as a proxy for the true answer, we find that our algorithm was generally within 2\% of the correct value which is acceptable performance. Interestingly, our algorithm performed worst for ``in-the-money'' options (such as when $Y_0 = 1.8, X_0 \sim \mathcal{N}(0.2, 0.05^2)$), i.e.\ when it is optimal to stop early. As expected, our method produced an underestimate of true $V$ since the computed stopping rule is necessarily sub-optimal. We found that the distribution of the computed $\tau^*$ was quite uniform on $\{\Delta t, 2\Delta t, \ldots, (M-1)\Delta t\}$ showing that this was a nontrivial stopping problem. For comparison, Table \ref{table:ex1} also lists the European option price assuming that early exercise is no longer possible. This column shows that our algorithm captured about 85-90\% of the \emph{time value of money}, i.e.\ the extra benefit due to early stopping.

\begin{figure}
 \centering{\includegraphics[height=2.4in,width=3in,clip]{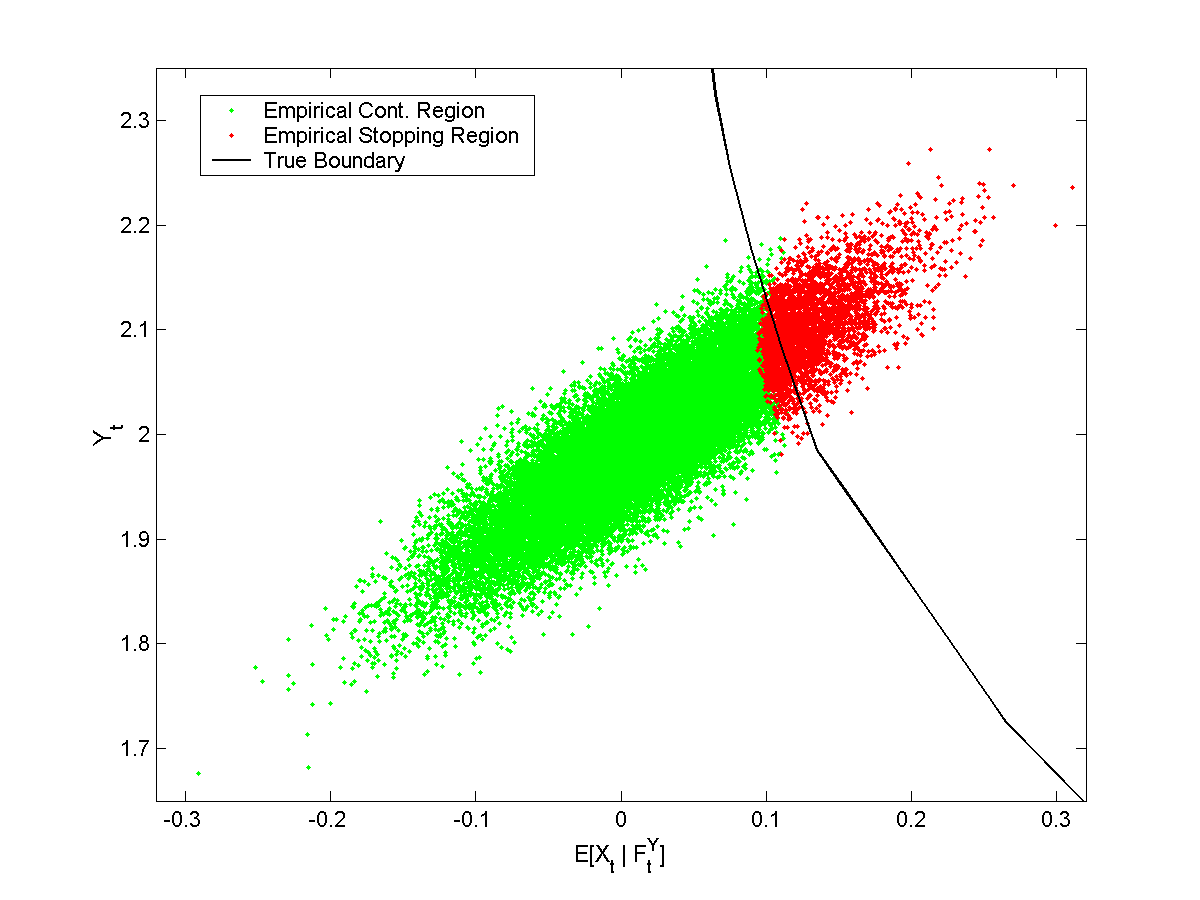}} \caption{Comparison of
optimal stopping regions for the pde and Monte Carlo solvers. The solid line shows the optimal stopping boundary as a function of $m_t$; the color-coded points show the values $z^n_j(t)$, $j=1,\ldots,N$ projected onto $\E[ X_t|\Fy_t] = (\rho_t x) \cdot (\rho_t 1)^{-1}$. Here $X_0 \sim \mathcal{N}(0, 0.05^2)$, $Y_0=2$, $N=30\,000$ and $t=0.5$. \label{fig:example1}}
\end{figure}

 To further illustrate the structure of the solution, Figure \ref{fig:example1} compares the optimal stopping regions computed by each algorithm at a fixed time point $t=0.5$. Note that since the value function $V$ is typically not very sensitive to the choice of a stopping rule, direct comparison of optimal stopping regions is more relevant (and more important for a practicing controller).  As we can see, an excellent fit was obtained through our non-parametric method. Figure \ref{fig:example1} also reveals that both $\{ q_t > G_t \} \cap \{ \hq_t \le G_t \}$ and $\{ q_t \le G_t \} \cap \{ \hq_t > G_t \}$ were non-empty (in other words, sometimes our algorithm stopped too early; sometimes it stopped too late). Recall that the simulation solver works under $\tP$ and therefore the empirical distribution of $Y_t$ in Figure \ref{fig:example1} would be different from the actual realizations under $\PP$ that will be observed by the controller.

While the pde formulation \eqref{eq:kalman-filter-equations}-\eqref{eq:example1-pde} is certainly better for the basic example above, it is crucially limited in its applicability. For instance, 
\eqref{eq:kalman-filter-equations} assumes Gaussian initial condition; any other $\xi_0$ renders it invalid. Similarly, perturbations to the dynamics \eqref{eq:model-1} will at the very least require re-derivation of \eqref{eq:kalman-filter-equations}-\eqref{eq:example1-pde}, or more typically lead to the case where no finite-dimensional sufficient statistics of $X_t | \Fy_t$ exist. In stark contrast to such difficulties, the particle filter algorithm can be used without any modifications for any $\xi_0$, and would need only minor adjustments to accommodate a different version of \eqref{eq:model-1}. A simple illustration is shown in the last two rows of Table \ref{table:ex1} where we consider a uniform and a discrete initial distribution, respectively. Heuristically, $V$ should be increasing with respect to the kurtosis of $\xi_0$, as a more spread-out initial distribution of $X_t$ leads to more optionality. Hence, (as confirmed by Table \ref{table:ex1}), $V(0,\xi^1,y_0) < V(0,\xi^2,y_0) < V(0,\xi^3,y_0)$, where $\xi^1 = \mathcal{N}(0,0.05^2), \xi^2 = Uniform([-0.05\sqrt{3},0.05\sqrt{3}]), \xi^3=0.5(\delta_{-0.05}+\delta_{0.05})$ are three initial distributions of $X$ normalized to $\int_\R x\xi^i(dx) = 0$, $\int_\R x^2 \xi^i(dx) = 0.05^2$.

\begin{table}
$$\begin{array}{|c|c|c|c|c|} \hline
\qquad \xi_0  & y_0  &  \text{Simulation solver} & \text{pde solver} & \text{European option}\\
\hline \mathcal{N}(0, 0.05^2) & 2 & 0.1810 & 0.1853 &  0.1331\\ 
\mathcal{N}(-0.12, 0.05^2) & 2.24 & 0.2566 & 0.2661 & 0.2136 \\ 
\mathcal{N}(0.2, 0.05^2) & 1.8 & 0.1862 & 0.1904  & 0.1052 \\
\mathcal{N}(0, 0.1^2) & 2 & 0.1852 & 0.1919 & 0.1349 \\
\delta_0 & 2 & 0.1723 & 0.1832 & 0.1325 \\
Unif_{[-0.05\sqrt{3},0.05\sqrt{3}]} & 2 & 0.1827 & -- & 0.1347 \\
0.5 (\delta_{-0.05} + \delta_{0.05}) & 2 & 0.1853 & -- &  0.1332\\ \hline
\end{array}$$
%
%
%
\caption{Comparison of the Monte Carlo scheme of Section \ref{sec:algo-summary} versus the Bermudan pde solver for the stochastic drift example of Section \ref{sec:examples}. Standard error of the Monte Carlo solver was about $0.001$.\label{table:ex1}}
\end{table}

\section{American Option Pricing under Stochastic Volatility}\label{sec:stoch-vol}
Our method can also be applied to \emph{stochastic volatility} models. Such asset pricing
models are widely used in financial mathematics to represent stock dynamics and assume that the local volatility of the underlying stock is itself stochastic. While under continuous observations the local volatility is perfectly known through
the quadratic variation process, under discrete observations this leads to a partially observed model
similar to \eqref{eq:basic-problem}.

To be concrete, let $Y_t$ represent the $\log$-price of a stock at time $t$ under the
given (pricing) measure $\PP$, and let $X_t$ be the instantaneous volatility of $Y$ at time $t$. We postulate that $(X,Y)$ satisfy the following system of sde's (known as the Stein-Stein model), 
\begin{align}\label{eq:stoch-vol} \left\{ \begin{aligned}
dY_t & = (r  - \half X^2_t)\, dt + X_t \,dU_t,\\
dX_t & = \kappa( \bar{\sigma} - X_t) \, dt + \rho \alpha \, dU_t + \sqrt{1-\rho^2}
\alpha \, dW_t. \end{aligned}\right.\end{align}
The stock price $Y$ is only
observed at the discrete time instances $\tT = \{ \Delta t, 2\Delta t, \ldots \}$ with $\tFy_t =
\sigma( Y_0, Y_{\Delta t}, \ldots, Y_{\lfloor t/\Delta t\rfloor \Delta t} )$. The American
(Put) option pricing problem consists in finding the optimal $\tFy$-adapted  and $\tT$-valued
stopping time $\tau$ for
\begin{align}\label{eq:obj-2}
\sup_{\tau \in \tT}\E[ \e^{-r \tau} (K - \e^{Y_\tau})_+ ].
\end{align}

A variant of \eqref{eq:stoch-vol}-\eqref{eq:obj-2} was recently studied by Sellami et al.\
\cite{SellamiPhamRunggaldier}. More precisely, \cite{SellamiPhamRunggaldier}  considered the American option pricing model in a simplified discrete
setting where $(X_t)$ of the Stein-Stein model \eqref{eq:stoch-vol} was replaced with a corresponding 3-state Markov chain approximation. In a related vein, Viens et al.\ \cite{Viens03} considered the filtering and portfolio optimization problem where the second line of \eqref{eq:stoch-vol} was replaced with the Heston model
\begin{align}\label{eq:heston} d(X^2_t) & = \kappa( \bar{\sigma} - X^2_t) \, dt + \rho \alpha X_t dU_t + \sqrt{1-\rho^2} \alpha {X_t} dW_t.
\end{align}
In general, the problem of estimation of $X_t$ is well-known, see e.g.\
\cite{CvitanicLiptserRozovskii06,FreyRunggaldier01,StroudPolson}. Observe that while \eqref{eq:stoch-vol} is linear, the square-root dynamics in \eqref{eq:heston} are highly non-linear and no finite-dimensional sufficient statistics exist for  $\pi_t$ in the latter case.

%
%

In the presence of stochastic volatility, one may no longer use the reference probability measure $\tP$. Indeed, there is no way to obtain a Brownian motion from the observation process $Y$ whose increments are now tied with the values of the unobserved $X$. Accordingly, $\zeta$ is no longer defined and consequently we cannot use it as an importance weight during the particle branching step in \eqref{eq:a-evolution}.  

A way out of this difficulty is provided by Del Moral et al.\ \cite{DelMoralProtter01}. The idea is to propagate particles independently of observations and to compute a \emph{candidate observation}
for each propagated particle. The weights are then assigned by comparing the candidates with the actual observation. Let $\phi$ be a smooth bounded function with $\int_\R \phi(x) \, dx = 1$ and $\int_\R |x| \phi(x) \,dx < \infty$ (e.g.\ $\phi(x) = \exp(-x^2/2) \cdot (2\pi)^{-1/2}$).
The propagated particles and candidates are obtained by
\begin{align*} \left\{ \begin{aligned}
v^n_j( t ) & = v^n_j( m\Delta t) + \int_{m\Delta t}^t \kappa(\bar{\sigma} - v^n_j(s) )\, ds + \int_{m\Dt}^t
\rho \alpha \, dU^{(j)}_s + \int_{m\Dt}^t \sqrt{1-\rho^2}\alpha \, dW^{(j)}_s,  \\
Y^{(j)}(t) & = \int_{m\Dt}^t (r  - \half (v^n_j(t))^2 )\, dt + \int_{m\Dt}^t v^n_j(s) \,dU^{(j)}_s, \qquad\qquad m\Dt \le t \le (m+1)\Dt,
 \end{aligned}\right.
 \end{align*}
where $(U^{(j)},W^{(j)})_{j=1}^n$ are $n$ independent copies of bivariate Wiener processes.
The branching weights are then given by
\begin{align}\label{eq:branching-protter}
\bar{a}^n_j( (m+1)\Dt) = n^{1/3} \phi \left( n^{1/3} (Y^{(j)}_{(m+1)\Delta t} - Y_{(m+1)\Delta t}) \right).
\end{align}
Hence, particles whose candidates $Y^{(j)}_{(m+1)\Delta t}$ are close to the true observed $Y_{(m+1)\Delta t}$ get high weights, while those particles that produced poor candidates are likely to be killed off. The rest of the algorithm remains the same as in Section \ref{sec:particle}. As shown in \cite[Theorem 5.1]{DelMoralProtter01}, the resulting filter satisfies for any bounded payoff function $f \in \cC^0_b(\R^d)$
\begin{align}\label{eq:protter-error}
\E[ | \pi^n_t f - \pi_t f| ] \le \frac{C(t)}{n^{1/3}} \| f \|_{0,\infty}, \qquad\qquad\text{with}\quad C(t) = \mathcal{O}(\e^t).
\end{align}
Note that compared to \eqref{eq:pi-error}, the error  in \eqref{eq:protter-error} as a function of number of particles $n$ is worse. This is due to the higher re-sampling variance produced by the additional randomness in $Y^{(j)}$'s. 

\subsection{Numerical Example}
Recently \cite{SellamiPhamRunggaldier} considered the above model \eqref{eq:stoch-vol} with 
the parameter values 
 $$ \begin{array}{c|ccccccccc}
\text{Parameter} & Y_0 &  X_0 & K  & \kappa & \bar{\sigma} & \alpha & T & r & \rho \\ \hline
\text{Value} & 110 & 0.15 & 100 & 1 & 0.15 & 0.1 & 1 & 0.05 & 0 \\
\end{array}.$$

Plugging-in the above parameters and using the modification \eqref{eq:branching-protter}, we implemented our algorithm with $N=30,000$, $n=1000$. Since no other solver of \eqref{eq:obj-2} is available, as in \cite{SellamiPhamRunggaldier} we compare the Monte Carlo solver of the partially-observed problem to a pde solver for the fully observed case (in which case the Bermudan option price is easily computed using the quasi-variational formulation based directly on \eqref{eq:stoch-vol}). Table \ref{table:stoch-vol} shows the results as we vary the observation frequency $\Dt$. Since $\Dt$ is also the frequency of the stopping decisions, smaller $\Dt$ increases both the partially and fully observed value functions. Moreover, as $\Dt$ gets smaller, the information set becomes richer and the handicap of partial information vanishes.
 
In this example where the payoff $K-\exp(Y_t)$ is a function of the observable $Y$ only, our algorithm obtains excellent performance. Also, we see that partial information has apparently only a mild effect on potential earnings (difference of less than 1.5\% even if $Y$ is observed just five times). To give an idea of the corresponding time value of money, the European option price in this example was $1.570$. Comparison with the results obtained in \cite{SellamiPhamRunggaldier} (first two columns of Table \ref{table:stoch-vol}) is complicated because the latter paper immediately discretizes $X$ and constructs a three-state Markov chain $(\tilde{X}_t)$. This discrete version takes on the values $\tilde{X}_{m\Dt} \in \{ 0.1, 0.15, 0.2 \}$ and therefore does not exhibit the asymmetric behavior of very small $X_t$ realizations that dampen the volatility of $Y$ and drastically reduce Put profits. In contrast, our algorithm operates on the original continuous-state formulation in \eqref{eq:stoch-vol}. Consequently, as can be seen in Table \ref{table:stoch-vol}, the full observation prices of the two models are quite different.

\begin{table}[h]
\begin{tabular}{|c|c|c|c|c|}\hline
 & \multicolumn{2}{c|}{Discrete Model} & \multicolumn{2}{c|}{Continuous Model} \\ 
$\Dt$ & Full Obs. & Partial Obs. & Full Obs. & Partial Obs. \\ \hline
0.2 & 1.575 & 0.988 &  1.665 & 1.646\\
0.1 & 1.726 & 1.306 & 1.686 &  1.673\\
0.05 & 1.912 & 1.596 & 1.696 & 1.685 \\ \hline
\end{tabular}
\caption{Comparison of discrete and continuous models for \eqref{eq:stoch-vol} under full and partial observations. The first two columns are reproduced from \cite[Table 3]{SellamiPhamRunggaldier}.\label{table:stoch-vol}}
\end{table}

%
%
%



\section{Conclusion}\label{sec:conclusion}
In this paper we have presented a new numerical scheme to solve partially observable optimal stopping problems. Our method is entirely simulation-based and only requires the ability to simulate the state processes. Consequently, we believe it is more robust than other proposals in the existing literature.

While our analysis was stated in the most simple setting of multi-dimensional diffusions, it can be considerably extended. First, as explained in Section \ref{sec:stoch-vol}, our algorithm can be easily adjusted to take into account \emph{discrete} observations which is often the more realistic setup. Second, the assumption of diffusion state processes is not necessary from a numerical point of view; one may consider other cases such as models with jumps, or even discrete-time formulations given in terms of general transition semigroups. For an example using a particle filter to filter a stable L\'evy process $X$, see \cite{KouritzinSun05}.  Third, one may straightforwardly incorporate \emph{state} constraints on the unobserved factor $X$. For instance, some applications imply that $X_t \ge 0$ is an extra constraint on top of \eqref{eq:X-sde} (in other words the observable filtration is generated by $Y$ and $1_{\{X_t \ge 0\}}$). Such a restriction can be added by assigning zero weights to particles that violate state constraints so that they are not propagated during the next branching step. Finally, if one uses the modification \eqref{eq:branching-protter} from \cite{DelMoralProtter01}  then many other noise formulations can be chosen beyond \eqref{eq:Y-sde}.

\bibliography{ospi}
\bibliographystyle{abbrv}

\end{document}